# MARTINGALE TRANSFORMS GOODNESS-OF-FIT TESTS IN REGRESSION MODELS[1]

By Estate V. Khmaladze and Hira L. Koul

*Victoria University of Wellington and Michigan State University*

This paper discusses two goodness-of-fit testing problems. The first problem pertains to fitting an error distribution to an assumed nonlinear parametric regression model, while the second pertains to fitting a parametric regression model when the error distribution is unknown. For the first problem the paper contains tests based on a certain martingale type transform of residual empirical processes. The advantage of this transform is that the corresponding tests are asymptotically distribution free. For the second problem the proposed asymptotically distribution free tests are based on innovation martingale transforms. A Monte Carlo study shows that the simulated level of the proposed tests is close to the asymptotic level for moderate sample sizes.

**1. Introduction.** This paper is concerned with developing asymptotically distribution free tests for two testing problems. The first problem pertains to testing a goodness-of-fit hypothesis about the error distribution in a class of nonlinear regression models. The second problem pertains to fitting a regression model in the presence of the unknown error distribution. The tests are obtained via certain martingale transforms of some residual empirical processes for the first problem and partial sum residual empirical processes for the second problem.

To be more precise, let $\Theta$ be an open subset of the $q$-dimensional Euclidean space and let $\{\mu(\cdot, \vartheta); \vartheta \in \Theta\}$ be a parametric family of functions from $\mathbb{R}^p$ to $\mathbb{R}$. For a pair $(X, Y)$ of a $p$-dimensional random vector $X$ with distribution function (d.f.) $H$ and one-dimensional random variable (r.v.) $Y$ with finite expectation let

$$m(x) := E[Y|X = x], \qquad x \in \mathbb{R}^p,$$

Received March 2002; revised February 2003.

[1]Supported in part by ARC grant RMS2974 and NSF Grant DMS-00-71619.

denote the regression function of $Y$ on $X$. In the first problem of interest one assumes $m$ is a member of a parametric family $\{\mu(\cdot, \vartheta); \vartheta \in \Theta\}$ and one observes a sequence $\{(X_i, Y_i), 1 \leq i \leq n\}$ such that for some $\theta \in \Theta$, the errors

$$\varepsilon_i(\theta) = Y_i - \mu(X_i, \theta), \qquad 1 \leq i \leq n, \tag{1.1}$$

are independent, identically distributed (i.i.d.) r.v.'s with expected value 0. Let $F$ be a specified distribution function with mean 0 and finite Fisher information for location, that is, $F$ is absolutely continuous with a.e. derivtive $f'$ satisfying

$$0 < \int \left(\frac{f'}{f}\right)^2 dF < \infty. \tag{1.2}$$

The problem of interest is to test the hypothesis

$$H_0 \colon \text{the d.f. of } \varepsilon_1(\theta) \text{ is } F,$$

against a class of all sequences of local (contiguous) alternatives where the error d.f.'s $A_n$ are such that for some $a \in L_2(\mathbb{R}, F)$,

$$\left(\frac{dA_n}{dF}\right)^{1/2} = 1 + \frac{1}{2\sqrt{n}} a + r_n,$$

$$\int a \, dF = 0, \tag{1.3}$$

$$n \int r_n^2 \, dF = o(1).$$

Occasionally, we will also insist that $a$ satisfy the orthogonality assumption

$$\int a \frac{f'}{f} \, dF = 0. \tag{1.4}$$

In the second problem one is again given independent observations $\{(X_i, Y_i), 1 \leq i \leq n\}$, such that $Y_i - m(X_i)$ are i.i.d. according to some distribution, not necessarily known, and one wishes to test the hypothesis

$$\widetilde{H}_0 \colon m(\cdot) = \mu(\cdot, \theta), \qquad \text{for some } \theta \in \Theta. \tag{1.5}$$

The alternative to $\widetilde{H}_0$ of interest here consists of all those sequences of functions $m_n(x)$ which "locally" deviate from one of $\mu(x, \theta)$, that is, for some $\theta \in \Theta$ and for some function $\ell_\theta \in L_2(\mathbb{R}^p, H)$,

$$\ell_\theta \perp \dot{\mu}_\theta, \qquad m_n(x) = \mu(x, \theta) + \frac{1}{\sqrt{n}} \ell_\theta(x) + r_{n\theta}(x),$$

$$n \int r_{n\theta}^2(x) \, dH(x) \to 0, \tag{1.6}$$

while the errors $Y_i - m_n(X_i)$ are still i.i.d. Here $\dot{\mu}_\theta(x)$ is a vector of $L_2$-derivatives of $\mu(x, \theta)$ with respect to $\theta$, assumed to exist; see the assumption (2.4).

Both of these testing problems are historically almost as old as the subject of statistics itself. The tests based on various residual empirical processes for $H_0$ have been discussed in the literature repeatedly. For example, see Durbin (1973), Durbin, Knott and Taylor (1975), Loynes (1980), D'Agostino and Stephens (1986) and Koul (1992, 2002), among others. Several authors have addressed the problem of regression model fitting, that is, testing for $\tilde{H}_0$: see, for example, Cox, Koh, Wahba and Yandell (1988), Eubank and Hart (1992, 1993), Eubank and Spiegelman (1990), Härdle and Mammen (1993), Koul and Ni (2004), An and Cheng (1991), Stute (1997), Stute, González Manteiga and Presedo Quindimil (1998), Stute, Thies and Zhu (1998) and Stute and Zhu (2002), among others. The last five references propose tests based on a certain marked empirical or partial sum processes while the former cited references base tests on nonparametric regression estimators. See also the review paper of MacKinnon (1992) for tests based on the least square methodology and the monograph of Hart (1997) and references therein for numerous other tests of $\tilde{H}_0$ based on smoothing methods in the case $p = 1$.

However, it is well known that most of these tests are not asymptotically distribution free. This is true even for the chi-square type of tests with the exception of the modified chi-square statistic studied in Nikulin (1973) in the context of empirical processes. It is also well documented in the literature that chi-square type tests often have relatively low power against many alternatives of interest, see, for example, Moore (1986). Hence a larger supply of asymptotically distribution free (ADF) goodness-of-fit tests with relatively good power functions is needed.

The aim of this paper is to propose a large class of such tests. These will be the tests based on statistics of a certain ADF modification and extension [see, e.g., (5.3) and (5.4)] of the (weighted) empirical process of residuals

$$\widehat{W}_\gamma(y) := n^{-1/2} \sum_{i=1}^n \gamma(X_i)[\mathbb{I}\{Y_i - \mu(X_i, \hat{\theta}) \leq y\} - F(y)],$$

$$-\infty \leq y \leq \infty,$$

where $\gamma$ is a square integrable function with respect to $H$. The ADF versions of the Cramér–von Mises and the Kolmogorov–Smirnov tests will be particular cases of such tests. Write $\widehat{W}_1$ for $\widehat{W}_\gamma$ whenever $\gamma \equiv 1$—see Sections 3.2 and 5.

As far as the problem of estimation of $\theta$ is concerned, certain weighted residual empirical processes play an indispensable role [cf. Koul (1992, 1996)].

A part of the objective of the present paper is to clarify the role of these processes with regard to the above goodness-of-fit testing problem.

To begin with, we shall discuss the basic structure of the first problem from a geometric perspective. This perspective was explored in the context of empirical processes in Khmaladze (1979). We shall show that under $H_0$ the asymptotic distribution of $\widehat{W}_\gamma$, and its general function-parametric form $\xi_n(\gamma, \varphi; \hat{\vartheta})$ [see (2.2)], is equivalent to that of the projection of (function-parametric) Brownian motion parallel to the tensor product $\dot{\mu}_\theta \cdot (f'/f)$. Since a "projection" is typically "smaller" than the original process we can intuitively understand why, at least for alternatives (1.3), it will lead to increase in asymptotic power if we substitute an estimator $\hat{\theta}$ even in the problems where the true value of the parameter is known. The distribution of this projection depends not only on the family of regression functions $\{\mu(\cdot, \vartheta); \vartheta \in \Theta\}$ and $F$, but also on the estimator $\hat{\theta}$. Therefore, the limit distribution of any fixed statistic based on $\widehat{W}_\gamma$ or on $\xi_n(\gamma, \varphi; \hat{\vartheta})$ will be very much model-dependent. However, using this "projection" point of view, we shall show in Section 3.2 that the tests based on $\widehat{W}_\gamma$ corresponding to a certain nonconstant $\gamma$ may be useful, because they may have simpler asymptotic behavior, but at the cost of some loss of the asymptotic power, and the tests based on $\widehat{W}_1$, in general, will have higher asymptotic power.

But, as mentioned above, the asymptotic null distribution of $\widehat{W}_1$ is model dependent. Proposed martingale transforms of $\widehat{W}_\gamma(F^{-1})$ will be shown to converge in distribution to a standard Brownian motion on $[0, 1]$ under $H_0$, and hence tests based on these transforms will be ADF for testing $H_0$. It will also be shown that for any $\gamma$ this transform is one-to-one and therefore there is no loss of the asymptotic power associated with it.

The paper also provides ADF tests for the problem of testing

$$H_\sigma: \text{ the d.f. of } \varepsilon_1(\theta) \text{ is } F(y/\sigma), \qquad \forall y \in \mathbb{R}, \text{ for some } \sigma > 0.$$

In the univariate design case, ADF tests for $\widetilde{H}_0$ based on certain partial sum processes and using ideas of Khmaladze (1981) have been discussed by Stute, Thies and Zhu (1998). An extension of this methodology to the general case of a higher dimensional design is far from trivial. The second important goal of this paper is to provide this extension. Here too we first discuss this problem from a general geometric perspective. It turns out that the weighted partial sum processes that are natural to this problem are

$$\xi_n(B; \hat{\theta}) := n^{-1/2} \sum \mathbb{I}\{X_i \in B\} \varphi(Y_i - \mu(X_i, \hat{\theta})),$$

for a fixed real valued function $\varphi$ with $E\varphi^2(\varepsilon)$ finite, where $B$ is a Borel measurable set in $\mathbb{R}^p$. Tests based on these processes and the innovation martingale transform ideas of Khmaladze (1993) [see, e.g., (6.4)] are shown to be ADF, that is, their asymptotic null distribution is free of the model $\mu(\cdot, \theta)$ and the error distribution, but depends on the design distribution in the case $p > 1$. These tests include those proposed in Stute, Thies and Zhu (1998), where $p = 1$, $\varphi(y) \equiv y$, $B = (-\infty, x]$, $x \in \mathbb{R}$.

We mention that recently Stute and Zhu (2002) used the innovation approach of Khmaladze (1981) to derive ADF tests in a special case of the higher dimension design where the design vector appears in the null parametric regression function only in a linear form, for example, as in generalized linear models, and where the sets $B$ in $\xi_n(B; \hat{\theta})$ are taken to be half spaces. This again reduces the technical nature of the problem to the univariate case.

In another recent paper Koenker and Xiao (2002) studied tests based on the transformations of a different process—regression quantile process to test the hypothesis that the effect of the covariate vector $X$ on the location and/or on the location-scale of the conditional quantiles of $Y$, given $X$, is linear in $X$. They then used the Khmaladze approach to make these tests ADF. Based on several Monte Carlo experiments, Koenker and Xiao (2001) report that their tests have accurate size and respectable power.

The paper is organized as follows. Section 2 introduces some basic processes that are used to construct tests of the above hypotheses. It also discusses some asymptotics under $H_0$ of these processes. Section 3 discusses some geometric implications of the asymptotics of Section 2, while Section 4 gives the martingale transforms of these processes whose asymptotic distribution under $H_0$ is known and free from $F$. Section 5 contains some computational formulas of these transformed processes. It also provides analogues of these ADF tests for nonrandom designs and when the underlying observations form a stationary autoregressive process. Section 6 contains the ADF processes for testing $\widetilde{H}_0$. Section 7 contains some simulation results to show how well the asymptotic level approximates the finite sample level for the proposed ADF tests. It is observed that even for the sample size 40, this approximation is quite good for the chosen simulation study. See Section 7 for details.

## 2. Function-parametric regression processes with estimated parameter.

2.1. *Function-parametric regression process.* Consider a regression process as is defined in Stute (1997):

$$\xi_n(B, y, \vartheta) := n^{-1/2} \sum_{i=1}^{n} \mathbb{I}\{X_i \in B\}[\mathbb{I}\{\varepsilon_i(\vartheta) \leq y\} - F(y)],$$

(2.1)
$$-\infty \leq y \leq \infty, \vartheta \in \Theta,$$

where $B$ is a Borel measurable set in the $p$-dimensional Borel space $(\mathbb{R}^p, \mathcal{B}(\mathbb{R}^p))$ and

$$\varepsilon_i(\vartheta) := Y_i - \mu(X_i, \vartheta), \qquad 1 \leq i \leq n.$$

We will use also notation $\mathbb{I}_B(X_i)$ for the indicator function $\mathbb{I}\{X_i \in B\}$ interchangeably. It is natural to consider an extension of the above process where the indicator weights are replaced by some weight function $\gamma(X_i)$. The function $\gamma$ may be scalar- or vector-valued. The weak convergence of such processes in the $y$ variable and for a fixed $\gamma$ has been developed in Koul ([1992](), [1996]()) and Koul and Ossiander ([1994]()).

It is not any less natural to consider an extension of these weighted empiricals to those processes where the second indicator involving the error random variable $\varepsilon_i(\vartheta)$ in (2.1) is also replaced by a function. Consider, therefore, a function-parametric version of (2.1) indexed by a pair of functions $(\gamma, \varphi)$:

$$
\begin{aligned}
(2.2) \quad \xi_n(\gamma, \varphi; \vartheta) &:= \int_{\mathbb{R}^{p+1}} \gamma(x)\varphi(y)\xi_n(dx, dy; \vartheta) \\
&= n^{-1/2} \sum_{i=1}^{n} \gamma(X_i) \Big[ \varphi(\varepsilon_i(\vartheta)) - \int \varphi(y)\, dF(y) \Big].
\end{aligned}
$$

We shall choose $\gamma \in L_2(\mathbb{R}^p, H)$ and $\varphi \in L_2(\mathbb{R}, F)$. In this way one can say that $\xi_n$ is defined for the function $\alpha(x, y) = \gamma(x)\varphi(y)$, which is an element of $\mathbb{L} := L_2(\mathbb{R}^{p+1}, H \times F)$. For a general $\alpha \in \mathbb{L}$ we certainly have

$$
\begin{aligned}
\xi_n(\alpha; \vartheta) &:= \int_{\mathbb{R}^{p+1}} \alpha(x, y)\xi_n(dx, dy; \vartheta) \\
&= n^{-1/2} \sum_{i=1}^{n} (\alpha(X_i, \varepsilon_i(\vartheta)) - E[\alpha(X_i, \varepsilon_i(\vartheta))|X_i]).
\end{aligned}
$$

We will realize, however, that it is sufficient and natural for our present purpose to restrict $\alpha$ to be of the above product type. In the sequel, for any functional $\mathcal{S}$ on $\mathbb{L}$ we will use the notation $\mathcal{S}(\alpha)$ or $\mathcal{S}(\gamma, \varphi)$ interchangeably, whenever $\alpha = \gamma \cdot \varphi$.

The processes defined at (2.1) and (2.2) are obviously closely related: (2.1) represents a regression process as a random measure on $\mathbb{R}^{p+1}$ while (2.2) represents it as an integral from this random measure. Also, (2.2) defines a linear functional on $\mathbb{L}$.

The function-parametric version (2.2) will help to visualize in a natural way the geometric picture of what is involved when we estimate parameters and show why and when we need "martingale transformations" (Sections 4 and 6) to obtain asymptotically distribution free tests.

2.2. *Asymptotic increments of $\xi_n$ with respect to parameter.* Since $\theta$ is unknown, in order to base tests of $H_0$ on the process $\xi_n$ we will need to replace it by an estimator $\hat{\theta}$ in this process. This estimator will be typically assumed to be $n^{1/2}$-consistent, that is,

$$
(2.3) \qquad\qquad \|\hat{\theta} - \theta\| = O_p(n^{-1/2}).
$$

There is thus a need to understand the behavior of $\xi_n(\alpha; \theta + n^{-1/2}v)$ as a process in $v \in \mathbb{R}^q, \|v\| \le k < \infty$. The first thing certainly is to consider the Taylor expansion of this function in $v$.

To do this assume the following $L_2$-differentiability condition of the regression function $\mu(x, \vartheta)$ with respect to $\vartheta$: there exists a $q \times 1$ vector $\dot{\mu}_\theta$ of functions from $\mathbb{R}^p \times \Theta$ to $\mathbb{R}^q$, such that

$$\mu(x, \vartheta) - \mu(x, \theta) = \dot{\mu}_\theta^T(x)(\vartheta - \theta) + \rho_\mu(x; \vartheta, \theta),$$

$$0 < \int \dot{\mu}_\theta^T(x)\dot{\mu}_\theta(x) \, dH(x) < \infty,$$

(2.4)

$$C_\theta := \int \dot{\mu}_\theta(x)\dot{\mu}_\theta^T(x) \, dH(x) \qquad \text{is positive definite,}$$

$$\int \sup_{\|\vartheta - \theta\| \le \epsilon} \rho_\mu^2(x; \vartheta, \theta) \, dH(x) = o(\epsilon^2), \qquad \text{as } \epsilon \to 0.$$

Here, and in the sequel, for any Euclidean vector $v$, $v^T$ denotes its transpose.

Now, if additionally $\varphi$ is differentiable with derivative $\varphi' \in L_2(\mathbb{R}, F)$ satisfying

$$(2.5) \qquad \lim_{\epsilon \to 0} \int \sup_{0 < \Delta < \epsilon} |\varphi'(y - \Delta) - \varphi'(y)|^2 \, dF(y) = 0,$$

then, with $\alpha(x, y) \equiv \gamma(x)\varphi(y)$, we have the following proposition.

PROPOSITION 2.1. *Under assumptions* (2.4) *and* (2.5), *the following holds for every* $0 < k < \infty$.

(i) *For any* $\gamma \in L_2(\mathbb{R}^p, H)$

$$\sup_{\|v\| \le k} |\xi_n(\alpha; \theta + n^{-1/2}v) - \xi_n(\alpha; \theta)$$

$$- E\gamma(X_i)\dot{\mu}_\theta^T(X)E\varphi'(\varepsilon)v| = o_p(1).$$

(ii) *For* $\gamma = \eta \mathbb{I}_B$, $B \in \mathcal{B}(\mathbb{R}^p)$ *and a fixed* $\eta \in L_2(\mathbb{R}^p, H)$,

$$\sup_{B \in \mathcal{B}, \|v\| \le k} \left| \xi_n(\alpha; \theta + n^{-1/2}v) - \xi_n(\alpha; \theta) \right.$$

$$\left. - n^{-1} \sum_{i=1}^n \eta_i \mathbb{I}_B(X_i)\dot{\mu}_\theta^T(X_i)\varphi'(\varepsilon_i)v \right| = o_p(1).$$

*Hence, under* (2.3) *one obtains*

$$(2.6) \quad \xi_n(\alpha; \hat{\theta}) = \xi_n(\alpha; \theta) - n^{-1} \sum_{i=1}^n \eta_i \mathbb{I}_B(X_i)\dot{\mu}_\theta^T(X_i)\varphi'(\varepsilon_i)n^{1/2}(\hat{\theta} - \theta) + \rho_n(B),$$

*where* $\rho_n(B)$ *is a sequence of stochastic processes indexed by* $B \in \mathcal{B}$, *tending to zero uniformly in* $B \in \mathcal{B}$ *in probability.*

The representation in (i) or in (2.6) will be very convenient and appropriate when dealing with the fitting of a regression model in Section 6. But for testing $H_0$ pertaining to the error distribution, as we will see in the next section, the differentiability of $\varphi$ is restrictive. We may wish, for example, to choose $\varphi$ to be an indicator function as in (2.1). Thus it is desirable to obtain an analog of the above proposition for as general a $\varphi$ as possible.

Towards this goal, let $\Phi$ denote the linear span of a class of nondecreasing real valued functions $\varphi(y)$, $y \in \mathbb{R}$, such that

$$\int \varphi^2(y)\, dF(y) < \infty,$$

(2.7)          $$\left( \int [\varphi(y-t) - \varphi(y-s)]^2\, dF(y) \right)^{1/2} \leq \nu(|t-s|),$$

$$-\epsilon \leq s, t \leq \epsilon,$$

for some $\epsilon > 0$ and for some continuous function $\nu$ from $[0, \infty)$ to $[0, \infty)$, with $\nu(0) = 0$, $\int_0^\epsilon \log \nu^{-1}(t)\, dt < \infty$. This is a wide class of functions and will be a source of our $\varphi$ in what follows.

For any two functions $\alpha$, $\beta \in \mathbb{L}$, let

$$\langle \alpha, \beta \rangle := \int_{\mathbb{R}^{p+1}} \alpha(x,y)\beta(x,y)\, dH(x)\, dF(y).$$

Note that if $\alpha$ or both $\alpha, \beta$ are vector functions, then $\langle \alpha, \beta \rangle$ or $\langle \alpha, \beta^T \rangle$ is a vector or a matrix of coordinate-wise inner products. Let $\|\alpha\| := \langle \alpha^T, \alpha \rangle^{1/2}$ for a vector function $\alpha$. Finally, let

$$\psi_f(y) := \frac{f'(y)}{f(y)},$$

$$m_\theta(x,y) := \dot{\mu}_\theta(x)\psi_f(y), \qquad x \in \mathbb{R}^p,\ y \in \mathbb{R}.$$

Note that

$$\langle \dot{\mu}_\theta, \dot{\mu}_\theta^T \rangle = C_\theta,$$

$$\langle m_\theta, m_\theta^T \rangle = C_\theta \|\psi_f\|^2.$$

We are ready to state

PROPOSITION 2.2.  *Suppose that* (1.2) *and* (2.4) *hold. Then for* $\alpha(x,y) = \gamma(x)\varphi(y)$ *with* $\gamma \in L_2(\mathbb{R}^p, H)$, $\varphi \in \Phi$,

(2.8)          $$\xi_n(\alpha; \hat{\theta}) = \xi_n(\alpha; \theta) + \langle \alpha, m_\theta^T \rangle n^{1/2}(\hat{\theta} - \theta) + o_p(1).$$

To appreciate some implications of (2.8) we need to consider those estimators of $\theta$ that admit an asymptotic linear representation. For the purpose of the present paper it would be enough to assume this. However, for completeness of the presentation we give a relatively broad set of sufficient conditions under which a class of M-estimators is asymptotically linear. Let $\{\eta_\vartheta, \vartheta \in \Theta\}$ be a family of $q$-dimensional functions on $\mathbb{R}^p$ with coordinates in $L_2(\mathbb{R}^p, H)$. Let $\beta_\vartheta := \eta_\vartheta \cdot \varphi, \ \varphi \in \Phi$. Define an M-estimator $\tilde{\theta}$ to be a solution of the equation

$$(2.9) \qquad\qquad \xi_n(\beta_\vartheta; \vartheta) = 0.$$

The following proposition gives a set of sufficient conditions for this estimator to be asymptotically linear.

PROPOSITION 2.3. *Suppose that* (1.2) *and* (2.4) *hold. In addition, suppose* $\varphi \in \Phi$ *and* $\{\eta_\vartheta, \vartheta \in \Theta\}$ *are such that*

$$(2.10) \qquad \int_{\mathbb{R}^p} \sup_{\|\vartheta - \theta\| < \epsilon} \|\eta_\vartheta - \eta_\theta\|^2 \, dH = o(1), \qquad as \ \epsilon \to 0,$$

*and the matrix* $\langle \beta_\theta, m_\theta^T \rangle$ *is nonsingular. Then* $\tilde{\theta}$ *defined at* (2.9) *satisfies*

$$(2.11) \qquad n^{1/2}(\tilde{\theta} - \theta) = -\langle \beta_\theta, m_\theta^T \rangle^{-1} \xi_n(\beta_\theta; \theta) + o_p(1).$$

*In particular, if* $\{\dot{\mu}_\vartheta; \vartheta \in \Theta\}$ *satisfies* (2.10), *then the solution* $\hat{\theta}$ *of the likelihood equation*

$$(2.12) \qquad\qquad \xi_n(m_\vartheta; \vartheta) = 0$$

*has the asymptotic linear representation*

$$(2.13) \qquad n^{1/2}(\hat{\theta} - \theta) = -\langle m_\theta, m_\theta^T \rangle^{-1} \xi_n(m_\theta; \theta) + o_p(1).$$

From now on $\hat{\theta}$ will stand for the solution of (2.12), and we shall use the abbreviated notation $\xi_n(\alpha) = \xi_n(\alpha; \theta)$, $\hat{\xi}_n(\alpha) = \xi_n(\alpha; \hat{\theta})$ and $\tilde{\xi}_n(\alpha) = \xi_n(\alpha; \tilde{\theta})$. Combining (2.8) with (2.11) and (2.13), we see that the leading term of $\hat{\xi}_n$ and of $\tilde{\xi}_n$, in general, can be represented as the linear transformation of $\xi_n$:

$$(2.14) \qquad \hat{\xi}_n(\alpha) = \xi_n(\alpha) - \langle \alpha, m_\theta^T \rangle \langle m_\theta, m_\theta^T \rangle^{-1} \xi_n(m_\theta) + o_p(1),$$

$$(2.15) \qquad \tilde{\xi}_n(\alpha) = \xi_n(\alpha) - \langle \alpha, m_\theta^T \rangle \langle \beta_\theta, m_\theta^T \rangle^{-1} \xi_n(\beta_\theta) + o_p(1).$$

These linear transformations have a remarkably simple and convenient structure as is described in Section 2.3.

2.3. *Processes $\hat{\xi}_n$ and $\tilde{\xi}_n$ as projections.* Let us use the notation 1 for the function in $y$ identically equal to 1, so that, for example, $\langle \varphi, 1 \rangle = \int \varphi(y)\, dF(y)$ and let $\varphi^1 = \varphi - \langle \varphi, 1 \rangle$, and for $\alpha = \gamma \cdot \varphi$ let $\alpha^1 = \gamma \cdot \varphi^1$. It is obvious that $\xi_n(\alpha) = \xi_n(\alpha^1)$.

For $\alpha \in \mathbb{L}$ and a vector-valued function $\beta$, with coordinates in $\mathbb{L}$, such that the matrix $\langle \beta, m_\theta^T \rangle$ is nonsingular (we require this for simplicity, although it is not necessary), let

$$(2.16) \qquad \Pi\alpha = \alpha - \langle \alpha, m_\theta^T \rangle \langle m_\theta, m_\theta^T \rangle^{-1} m_\theta,$$

$$(2.17) \qquad \Pi_\beta \alpha = \alpha - \langle \alpha, m_\theta^T \rangle \langle \beta, m_\theta^T \rangle^{-1} \beta.$$

PROPOSITION 2.4. (i) *The linear transformation $\alpha \mapsto \alpha^1$ is an orthogonal projection in $\mathbb{L}$ parallel to functions which are constant in $y$.*

(ii) *The linear transformation $\Pi_{\beta_\theta}$ (and therefore $\Pi$) is a projection. It projects parallel to $\beta_\theta$ on a subspace of functions orthogonal to $m_\theta$. In particular $\Pi$ is an orthogonal projection parallel to $m_\theta$.*

(iii) *Adjoint projectors $\Pi_{\beta_\theta}^*$ (and therefore $\Pi^*$) project parallel to $m_\theta$. For any two vector functions $\beta, \lambda$,*

$$(2.18) \qquad \Pi_{\bar\beta}^* \Pi_{\bar\lambda}^* = \Pi_{\bar\beta}^*.$$

We can therefore say that under the regularity conditions that guarantee the validity of the expansions at (2.8), (2.11) and (2.13), the substitution of the M-estimator $\bar{\theta}$ in $\xi_n(\alpha; \theta)$ for $\theta$ is asymptotically equivalent to projecting $\xi_n(\alpha; \theta)$ parallel to the linear functional $m_\theta$ generated by $\dot\mu_\theta$ and $\psi_f$. Similarly, the substitution of the MLE $\hat{\theta}$ in $\xi_n(\alpha; \theta)$ for $\theta$ is asymptotically equivalent to projecting $\xi_n(\alpha; \theta)$ orthogonal to $m_\theta$. Moreover, the property (2.18) shows that the leading terms of $\xi_n(\alpha; \tilde{\theta}_1)$ and $\xi_n(\alpha; \tilde{\theta}_2)$, for any two estimators $\tilde{\theta}_1, \tilde{\theta}_2$ admitting the asymptotic linear representation (2.11), are in one-to-one correspondence with each other. Even though one of the estimators may be asymptotically more efficient than the other, (2.18) shows that the stocks of test statistics based on each of these processes are asymptotically the same. Therefore the inference based on either $\xi_n(\alpha; \tilde{\theta}_1)$ or $\xi_n(\alpha; \tilde{\theta}_2)$ will be asymptotically indistinguishable.

We end this section by outlining the proofs for Propositions 2.2 and 2.4. Throughout, $\varepsilon_i$ stands for $\varepsilon_i(\theta)$, $1 \le i \le n$.

2.4. *Some proofs.*

SKETCH OF THE PROOF OF PROPOSITION 2.1. We shall sketch details only for part (ii), while those for part (i) are similar and simpler. Let $\Delta_i(v) =$

$\mu(X_i, \theta + n^{-1/2}v) - \mu(X_i, \theta)$. Rewrite

$$\xi_n(\alpha; \theta + n^{-1/2}v) - \xi_n(\alpha; \theta) - n^{-1} \sum_{i=1}^{n} \eta_i \mathbb{I}_B(X_i) \dot{\mu}_\theta^T(X_i) \varphi'(\varepsilon_i) v$$

$$= n^{-1} \sum_{i=1}^{n} \eta_i \mathbb{I}_B(X_i) [\varphi(\varepsilon_i + \Delta_i(v)) - \varphi(\varepsilon_i) - \Delta_i \varphi'(\varepsilon_i)]$$

$$+ n^{-1/2} \sum_{i=1}^{n} \eta_i \mathbb{I}_B(X_i) [\Delta_i(v) - \dot{\mu}_\theta^T(X_i) \, n^{-1/2}v] \varphi'(\varepsilon_i).$$

The condition (2.4) implies that for every $\epsilon > 0$, $\exists N_\epsilon < \infty$ such that with probability at least $1 - \epsilon$ the following holds for all $n > N_\epsilon$:

$$E\left\{ \sup_{\|v\| \leq k} \sum_{i=1}^{n} |\Delta_i(v) - n^{-1/2} \dot{\mu}_\theta^T(X_i) v|^2 \right\} \to 0, \qquad \sup_{1 \leq i \leq n; \|v\| \leq k} |\Delta_i(v)| = o_p(1).$$

This fact and (2.5) imply the conclusion (ii) in a routine fashion. $\square$

Before proving the next proposition, we recall from Hájek (1972) that (1.2) implies the mean-square differentiability of $f^{1/2}$:

$$\frac{f^{1/2}(y + \delta) - f^{1/2}(y)}{f^{1/2}(y)} = \frac{1}{2} \frac{f'}{f}(y)\delta + \rho_f(y; \delta),$$

$$\int \rho_f^2(y; \delta) \, dF(y) = o(\delta^2), \qquad \delta \to 0.$$

This fact is used implicitly in the following proof and throughout the discussion in the paper without mentioning it explicitly.

PROOF OF PROPOSITION 2.2. Recall $\alpha(x, y) = \gamma(x)\varphi(y)$. Rewrite $\xi_n = \xi_{no} + \xi_n^*$, where

$$\xi_{no}(\alpha; \vartheta) := n^{-1/2} \sum_{i=1}^{n} \gamma(X_i) \Big[ \varphi(\varepsilon_i(\vartheta)) - E_\theta[\varphi(\varepsilon_i(\vartheta))|X_i] \Big],$$

$$\xi_n^*(\alpha; \vartheta) := n^{-1/2} \sum_{i=1}^{n} \gamma(X_i) \Big[ E_\theta[\varphi(\varepsilon_i(\vartheta))|X_i] - E_\theta[\varphi(\varepsilon_i(\theta))] \Big].$$

Note that $\xi_{no}(\alpha; \theta) = \xi_n(\alpha; \theta)$.

To prove Proposition 2.2 it thus suffices to show that for every $0 < k < \infty$,

(2.19)
$$\sup_{\|v\| \leq k} |\xi_{no}(\alpha; \theta + n^{-1/2}v) - \xi_{no}(\alpha; \theta)| = o_p(1),$$

(2.20)
$$\sup_{\|v\| \leq k} |\xi_n^*(\alpha; \theta + n^{-1/2}v) - m^T(\alpha; \theta) n^{-1/2} u| = o_p(1).$$

But (2.19) will follow from the equicontinuity condition of the process $\xi_{no}(\alpha; \cdot)$:

$$\sup_{\|\vartheta - \theta\| \leq \epsilon} |\xi_{no}(\alpha; \vartheta) - \xi_{no}(\alpha; \theta)| = o_p(1),$$

as $n \to \infty$ and $\epsilon \to 0$. This in turn follows from the argument below.

A $\varphi \in \Phi$ may be written as $\varphi = \varphi_1 - \varphi_2$, where nondecreasing $\varphi_1$, $\varphi_2$ both satisfy (2.7). Let $I_i := \text{sign}(\gamma(X_i))$, $i = 1, \ldots, n$. Then for any $\delta > 0$ and for all $i = 1, \ldots, n$,

$$\gamma(X_i)[\varphi_1(Y_i - \Delta - \delta I_i) - \varphi_2(Y_i - \Delta + \delta I_i)]$$
$$\leq \gamma(X_i)\varphi(Y_i - \Delta)$$
$$\leq \gamma(X_i)[\varphi_1(Y_i - \Delta + \delta I_i) - \varphi_2(Y_i - \Delta - \delta I_i)].$$

The expected value of the square of the above upper and lower bounds is bounded from above by

$$\int \gamma^2 \, dH \, 2\nu^2(2\delta).$$

Therefore the bracketing entropy (log of covering number) does not exceed

$$\log \nu^{-1} \left( t \Big/ \left[ 2 \int \gamma^2 \, dH \right]^{1/2} \right),$$

and hence is integrable by the definition of $\nu$. Therefore, by a result in van der Vaart and Wellner (1996, Sections 2.5.2, 2.7), (2.19) follows.

To prove (2.20), let, as above, $\Delta_i(v) = \mu(X_i, \theta + n^{-1/2}v) - \mu(X_i, \theta)$. Then one has

$$\xi_n^*(\alpha; \theta + n^{-1/2}v)$$
$$= n^{-1/2} \sum_{i=1}^n \gamma(X_i) \int \varphi(y)[f(y + \Delta_i(v)) - f(y)] \, dy$$
$$= n^{-1} \sum_{i=1}^n \gamma(X_i)\dot{\mu}^T(X_i) \int \varphi(y)\psi_f(y) \, dF(y)v + \rho_n(v)$$
$$= \langle \alpha, m_\theta^T \rangle v + \rho_n^*(v),$$

where under the assumed conditions and using an argument similar to one used, for example, in Hájek and Šidák (1967) one can show that $\sup_{\|v\| \leq k} |\rho_n^*(v)| = o_p(1)$. $\qquad \square$

PROOF OF PROPOSITION 2.4. Let us prove part (iii) only. We need to show that $\Pi_\beta^* \Pi_\lambda^* \mathcal{S}(\alpha) = \Pi_\beta^* \mathcal{S}(\alpha)$. We have

$$\Pi_\beta^* \Pi_\lambda^* \mathcal{S}(\alpha) = \Pi_\beta^* \mathcal{S}(\alpha) - \Pi_\beta^* \langle \alpha, m_\theta^T \rangle \langle \lambda, m_\theta^T \rangle^{-1} \mathcal{S}(\lambda).$$

But, by definition,

$$\Pi_\beta^* \langle \alpha, m_\theta^T \rangle = \langle \alpha, m_\theta^T \rangle - \langle \alpha, m_\theta^T \rangle \langle \beta, m_\theta^T \rangle^{-1} \langle \beta, m_\theta^T \rangle = 0.$$

Hence the last claim. It implies that $\Pi_\beta^* \Pi_\beta^* = \Pi_\beta^*$, that is, $\Pi_\beta^*$ is a projection. The remainder of the proof is obvious.   $\square$

### 3. Limiting process and asymptotic power.

3.1. *The limiting process.* Let $b(x, y), x \in \mathbb{R}^p, y \in \mathbb{R}$, be a Brownian motion with covariance function $H(x \wedge x')F(y \wedge y')$, where $x \wedge x'$ is the vector with coordinates $\min(x_i, x_i'), i = 1, \ldots, p$. In the discussion below all $\gamma$'s and $\varphi$'s are in $L_2(\mathbb{R}^p, H)$ and $L_2(\mathbb{R}, F)$, respectively, that is, $(\gamma, \varphi) \in \mathbb{L}$. Define, for $\alpha(x, y) = \gamma(x)\varphi(y)$, the function parametric Brownian motion

$$b(\alpha) := b(\gamma, \varphi) := \int_{\mathbb{R}^{p+1}} \gamma(x)\varphi(y)\, b(dx, dy).$$

Clearly the class $\{b(\alpha) \colon \alpha \in \mathbb{L}\}$ is a family of zero mean Gaussian random variables with the covariance given by

$$Eb(\alpha_1)b(\alpha_2) = \langle \alpha_1, \alpha_2 \rangle.$$

Let

$$\xi(\alpha) := b(\gamma, \varphi) - \langle \varphi, 1 \rangle b(\gamma, 1) = b(\gamma, \varphi^1) = b(\alpha^1).$$

The family $\{\xi(\alpha) \colon \alpha \in \mathbb{L}\}$ is also a family of zero mean Gaussian random variables with the covariance

$$E\xi(\alpha_1)\xi(\alpha_2) = \langle \gamma_1, \gamma_2 \rangle [\langle \varphi_1, \varphi_2 \rangle - \langle \varphi_1, 1 \rangle \langle \varphi_2, 1 \rangle] = \langle \alpha_1^1, \alpha_2^1 \rangle.$$

Thus, $\xi(\alpha)$ is a function parametric Kiefer process in $\alpha$ and simply a Brownian motion in $\alpha^1$. Finally, define

$$\hat\xi(\alpha) := \xi(\alpha) - \langle \alpha, m_\theta^T \rangle \langle m_\theta, m_\theta^T \rangle^{-1} \xi(m_\theta) = \Pi\xi(\alpha).$$

Since $\langle \psi_f, 1 \rangle = \int f'(y)\, dy = 0$, we have $\xi(m_\theta) = b(m_\theta)$. Hence, $\hat\xi$ can be rewritten as

$$(3.1) \quad \hat\xi(\alpha) = b(\alpha^1) - \langle \alpha^1, m_\theta^T \rangle \langle m_\theta, m_\theta^T \rangle^{-1} b(m_\theta) = \Pi b(\alpha^1) = b(\Pi\alpha^1).$$

It seems easier to use below the notation $\alpha_\perp$ for $\Pi\alpha^1$:

$$\alpha_\perp = \alpha^1 - \langle \alpha^1, m_\theta \rangle \langle m_\theta, m_\theta^T \rangle^{-1} m_\theta,$$

which is the part of $\alpha$ orthogonal to 1 and $m_\theta$.

Here and everywhere below we will consider only the case of orthogonal projectors, which asymptotically correspond to the substitution of the MLE. As our comment after Proposition 2.4 shows, we can do this without loss of generality.

In view of (2.14), the reason for introducing the processes $\xi$ and $\hat\xi$ is clear and is given by the following statement.

PROPOSITION 3.1. *Suppose that the conclusion* (2.14) *holds. Then the following holds for every* $\alpha \in \mathbb{L}$:

*Under* $H_0$

$$(3.2) \qquad \xi_n(\alpha) \overset{d}{\to} \xi(\alpha), \qquad \hat{\xi}_n(\alpha) \overset{d}{\to} \hat{\xi}(\alpha).$$

*Under the alternatives* (1.3)

$$(3.3) \qquad \begin{aligned} &\xi_n(\alpha) \overset{d}{\to} \xi(\alpha) + \langle \alpha, a \rangle, \\ &\hat{\xi}_n(\alpha) \overset{d}{\to} \hat{\xi}(\alpha) + \langle \alpha, a \rangle - \langle \alpha, m_\theta^T \rangle \langle m_\theta, m_\theta^T \rangle^{-1} \langle m_\theta, a \rangle. \end{aligned}$$

Because both $\xi_n$ and $\hat{\xi}_n$ are linear in $\alpha$, the above proposition is equivalent to the weak convergence of any finite-dimensional distributions of these processes. Hence the possible weak limits of these processes are uniquely determined.

From (3.2) and (3.3) we see that the asymptotic shift of $\hat{\xi}_n$ under the alternatives (1.3) and (1.4) is simply $\langle \alpha, a \rangle$, if $\alpha \perp m_\theta$, that is, if either $\gamma \perp \dot{\mu}_\theta$ or $\varphi \perp \psi_f$.

3.2. *The case of* $\gamma \perp \dot{\mu}_\theta$. In this case there exists an optimal choice of $\gamma$ which will maximize the asymptotic "signal to noise" ratio $\Delta$ of $\hat{\xi}_n(\gamma, \varphi)$ uniformly in $a$, that is, uniformly in alternatives (1.3), where

$$\Delta := \frac{|\langle \alpha, a \rangle|}{\|\alpha\|} = \frac{|\langle \gamma, 1 \rangle|}{\|\gamma\|} \frac{|\langle \varphi, a \rangle|}{\|\varphi\|}.$$

Here, too, we use the notation 1 for the function in $x$ identically equal to 1. Clearly, the $\gamma$ that maximizes $\Delta$, uniformly in $a$, is the $\gamma$ that maximizes the ratio

$$(3.4) \qquad \frac{|\langle \gamma, 1 \rangle|}{\|\gamma\|},$$

subject to the condition that $\gamma \perp \dot{\mu}_\theta$, and is given by

$$1_\perp := 1 - \int \dot{\mu}_\theta^T \, dH \, C_\theta^{-1} \dot{\mu}_\theta = 1 - \langle \dot{\mu}_\theta^T, 1 \rangle C_\theta^{-1} \dot{\mu}_\theta.$$

On the other hand, the $\gamma$ that maximizes $\Delta$ or (3.4) among all $\gamma \in L_2(\mathbb{R}^p, H)$ is 1. Then $1_\perp$ is simply part of the identity function 1 orthogonal to $\dot{\mu}_\theta$. It thus follows that $1_\perp(x) \equiv 0$ when $\mu(x, \vartheta)$ is linear in $\vartheta$ and has a nonzero intercept.

Now consider $\hat{\xi}(1_\perp, \varphi)$ as a process in $\varphi$, assuming that $\|1_\perp\| \neq 0$. Since $1_\perp \cdot \varphi$ is orthogonal to $\dot{\mu}_\theta$, from (3.1) we obtain

$$\hat{\xi}(1_\perp, \varphi) = \xi(1_\perp, \varphi) = b(1_\perp, \varphi^1).$$

It thus follows that $\hat{\xi}(1_\perp, \varphi)$ is a Brownian bridge in $\varphi$. If, for example, we choose $\varphi(y) = \varphi_t(y) = \mathbb{I}(y \leq F^{-1}(t))$, $0 \leq t \leq 1$, then along the family of functions $\{\varphi_t(\cdot), 0 \leq t \leq 1\}$, the process

$$u(t) := \hat{\xi}\left(\frac{1_\perp}{\|1_\perp\|}, \varphi_t\right)$$

is a standard Brownian bridge with $Eu(s)u(t) = s \wedge t - st$.

A prelimiting form of the process $u$ is

$$\hat{u}_n(t) = \hat{\xi}_n\left(\frac{1_{\perp,n}}{\|1_{\perp,n}\|_n}, \varphi_t\right)$$

$$= n^{-1/2} \sum_{i=1}^n \frac{1_{\perp,n}(X_i)}{\|1_{\perp,n}\|_n} [\mathbb{I}\{\varepsilon_i(\hat{\theta}) \leq F^{-1}(t)\} - t],$$

$$1_{\perp,n}(x) := [1 - \langle \dot{\mu}_{\hat{\theta}}^T, 1 \rangle_n C_{\hat{\theta},n}^{-1} \dot{\mu}_{\hat{\theta}}(x)], \qquad x \in \mathbb{R}^p,$$

where

$$C_{\hat{\theta},n} := \int \dot{\mu}_{\hat{\theta}} \dot{\mu}_{\hat{\theta}}^T \, dH_n,$$

$$\langle \dot{\mu}_{\hat{\theta}}^T, 1 \rangle_n = \int \dot{\mu}_{\hat{\theta}}^T \, dH_n,$$

$$\|1_{\perp,n}\|_n := (1 - \langle \dot{\mu}_{\hat{\theta}}^T, 1 \rangle_n C_{\hat{\theta},n}^{-1} \langle \dot{\mu}_{\hat{\theta}}, 1 \rangle_n)^{1/2}$$

and where $H_n$ is the empirical d.f. of the design variables $\{X_i, 1 \leq i \leq n\}$. One can verify, using, for example, the results from Koul (1996), that under the present setup, $\hat{u}_n$ converges weakly to a Brownian bridge. Hence, for instance, tests based on

$$\sup_t |\hat{u}_n(t)| \quad \text{or} \quad \int_0^1 |\hat{u}_n(t)|^2 \, dt$$

will have asymptotically the well-known Kolmogorov and Cramér–von Mises distributions, respectively.

Now, suppose that the design d.f. $H$ and the regression function $\mu_\theta$ are such that

$$(3.5) \qquad \langle \dot{\mu}_\theta, 1 \rangle = \int \dot{\mu}_\theta \, dH = 0.$$

Then $1_\perp \equiv 1$, and $\hat{u}_n \equiv \widehat{W}_1(F^{-1})$, the ordinary empirical process of the residuals whose weak convergence to Brownian bridge can also be derived from Koul (1996) under the present setup.

There is, however, a drawback in the choice of $\gamma \perp \dot{\mu}_\theta$: although, as we see, this choice of $\gamma$ makes the asymptotic behavior of $\hat{\xi}_n$ in $\varphi$ simple, the tests based on the process $\hat{\xi}_n(\|1_\perp\|^{-1}1_\perp, \varphi)$ will in general have some loss of asymptotic power. Consider for the moment the problem of testing $H_0$ vs. the alternative (1.3) for given $a$ when $\theta$ is known. Then the shift function that will appear in the asymptotic power for $\xi_n(\gamma, \varphi)$ is $|\langle\gamma, 1\rangle\langle\varphi, a\rangle|/\|\gamma\|\|\varphi\|$. This will attain its maximum in $\gamma$ when $\gamma \equiv 1$. However, for the process $\hat{\xi}_n(1_\perp, \varphi)$ the corresponding shift is uniformly smaller in absolute value:

$$\left|\frac{\langle 1_\perp, 1\rangle}{\|1_\perp\|}\langle\varphi, a\rangle\right| < |\langle\varphi, a\rangle|$$

and, in particular, the statistic $\hat{\xi}_n(1_\perp, a)$ will have smaller asymptotic power against the alternative $a$ than the statistic $\xi_n(1, a)$. The actual loss may be quite small, depending on the quantity

$$\|1_\perp\|^2 = 1 - \langle\dot\mu_\theta^T, 1\rangle C_\theta^{-1}\langle\dot\mu_\theta^T, 1\rangle,$$

and may actually equal 0, if (3.5) holds. But, in general, there is some loss.

We shall see in Section 6 that the choice of $\gamma \perp \dot\mu_\theta$ will become most natural when fitting a regression model. However, one should not think that the loss of power associated with this choice in testing the hypothesis $H_0$ is unavoidable due to the estimation of the nuisance parameters. On the contrary, estimation of the parameter may lead to an *increase of power* against "most" alternatives. We will see this better in the next section.

Finally, we remark that the geometric picture, similar to the one depicted by Propositions 2.4 and 3.1 and also in this and the next sections, was developed in the context of the parametric empirical processes in Khmaladze (1979). See also the monograph by Bickel, Klaassen, Ritov and Wellner (1998) describing the related geometry in connection with efficient and adaptive estimation in semiparametric models.

### 3.3. *The case of $\varphi \perp \psi_f$.*

This case is important for two reasons. The first is that in this case again $\hat\xi(\gamma, \varphi) = \xi(\gamma, \varphi)$, that is, the asymptotic behavior of the processes $\hat\xi_n(\alpha)$ and $\xi_n(\alpha)$ under $H_0$ is the same. The second is that if we assume that $a$ of (1.3) also satisfies (1.4), then there is in general a gain in the signal to noise ratio if we choose $\varphi$ orthogonal to $\psi_f$. Indeed, let $\varphi_\perp$ denote the part of $\varphi$ orthogonal to $\psi_f$ and 1. The signal to noise ratio for $\hat\xi_n(\gamma, \varphi_\perp)$ is *asymptotically larger* than that for $\xi_n(\gamma, \varphi)$, as is seen from the following elementary argument:

$$\frac{\langle\gamma, 1\rangle}{\|\gamma\|}\frac{\langle\varphi, a\rangle}{\|\varphi^1\|} = \frac{\langle\gamma, 1\rangle}{\|\gamma\|}\frac{\langle\varphi_\perp, a\rangle}{\|\varphi^1\|} \leq \frac{\langle\gamma, 1\rangle}{\|\gamma\|}\frac{\langle\varphi_\perp, a\rangle}{\|\varphi_\perp\|},$$

because $\|\varphi^1\| \geq \|\varphi_\perp\|$.

It is also obvious that the optimal choice of $\gamma$ that maximizes $\Delta$ uniformly in $a$ is $\gamma \equiv 1$. Therefore, consider the process

$$(3.6) \qquad\qquad \hat\xi(1, \varphi) = \xi(1, \varphi) = b(1, \varphi)$$

as a process in $\varphi$, for $\varphi$ satisfying $\varphi \perp \psi_f$ and $\varphi \perp 1$. From (3.6) it is clear that if we had a family of functions $\{\varphi_t, 0 \leq t \leq 1\}$ from $L_2(\mathbb{R}, F)$ such that

$$(3.7) \qquad \langle \varphi_t, 1 \rangle = \langle \varphi_t, \psi_f \rangle = 0,$$

$$(3.8) \qquad \langle \varphi_t, \varphi_t \rangle = t, \qquad 0 \leq t \leq 1,$$

$$(3.9) \qquad \langle \varphi_{t_2} - \varphi_{t_1}, \varphi_{t_1} \rangle = 0, \qquad t_2 \geq t_1,$$

then the process $\xi(1, \varphi_t)$, $0 \leq t \leq 1$, would be a Brownian motion in $0 \leq t \leq 1$. Hence, all tests based on

$$n^{-1/2} \sum_{i=1}^{n} \varphi_t(\varepsilon_i(\hat{\theta})), \qquad 0 \leq t \leq 1,$$

will be ADF.

It is straightforward to construct a family of functions satisfying (3.8) and (3.9). For example, take any function $\varphi$ from $L_2(\mathbb{R}, F)$ such that $L(y) := \int_{-\infty}^{y} \varphi^2 \, dF$ is a continuous distribution function on $\mathbb{R}$, and $\varphi^2 f > 0$, a.e. Then the family

$$(3.10) \qquad \begin{aligned} \varphi_t(y) &:= \varphi(y) \mathbb{I}\{y \leq L^{-1}(t)\}, \\ L^{-1}(t) &:= \inf\{y \in \mathbb{R} : L(y) \geq t\}, \qquad 0 \leq t \leq 1, \end{aligned}$$

satisfies these conditions. However, finding a family $\{\varphi_t, 0 \leq t \leq 1\}$ that satisfies (3.7) as well becomes far less straightforward. It is here we will exploit the "martingale transform" ideas of Khmaladze (1981, 1993).

## 4. A martingale transform.

Let $h(y) := (1, \psi_f(y))^T$ be an extended score function of the error distribution and set

$$\Gamma_t := \int_{z \geq y} h(z) h^T(z) \, dF(z) = \begin{pmatrix} 1 - F(y) & -f(y) \\ -f(y) & \int_y^\infty \psi_f^2(z) \, dF(z) \end{pmatrix},$$
$$t = F(y).$$

The matrix $\Gamma_t$ will be assumed to be nonsingular for every $0 \leq t < 1$. This, indeed, is true if and only if 1 and $\psi_f(y)$ are linearly independent on the set $y > c$ for all sufficiently large $c$. This, in turn, is true if $\psi_f$ is not a constant in the right tail of the support of $f$. Then the unique inverse $\Gamma_t^{-1}$ exists for every $0 \leq t < 1$. [The case when $\Gamma_t$ is not uniquely invertible does not create, however, much of a problem for the transformation (4.1), as is shown in Tsigroshvili (1998).]

Now, observe that the condition (3.7) above is equivalent to requiring that $\varphi$ be orthogonal to the vector $h$. For a function $\varphi \in L_2(\mathbb{R}, F)$, consider the transformation

$$(4.1) \qquad \mathcal{L}\varphi(y) := \varphi(y) - \int_{z \leq y} \varphi(z) h^T(z) \Gamma_{F(z)}^{-1} \, dF(z) \, h(y), \qquad y \in \mathbb{R}.$$

Let, for a $(\gamma, \varphi) \in \mathbb{L}$,

$$w(\alpha) := \hat{\xi}(\gamma, \mathcal{L}\varphi), \qquad \alpha = \gamma \cdot \varphi.$$

We have the following:

PROPOSITION 4.1. *Let* $\mathcal{H} := \{\varphi \in L_2(\mathbb{R}, F) \colon \langle \varphi, h \rangle = 0\}$. *The transformation* $\mathcal{L}$ *of* (4.1) *is a norm preserving transformation from* $L_2(\mathbb{R}, F)$ *to* $\mathcal{H}$:

$$\mathcal{L}\varphi \perp h, \qquad \|\mathcal{L}\varphi\| = \|\varphi\|.$$

*Consequently the process* $w(\alpha)$ *is a (function parametric) Brownian motion on* $\mathbb{L}$.

A consequence of this proposition is the following corollary:

COROLLARY 4.1. *Suppose* $\{\varphi_t, 0 \le t \le 1\}$ *is a family of functions satisfying the conditions* (3.8) *and* (3.9). *Then* $\{\mathcal{L}\varphi_t, 0 \le t \le 1\}$ *is a family of functions satisfying all three conditions* (3.7)–(3.9). *Consequently,* $\{\hat{\xi}(\gamma, \mathcal{L}\varphi_t), 0 \le t \le 1\}$, *for any fixed* $\gamma$ *with* $\|\gamma\| = 1$, *is a standard Brownian motion in* $t$.

Now, if $\{\hat{\xi}_n(\gamma, \mathcal{L}\varphi_t), 0 \le t \le 1\}$ converges weakly to $\{\hat{\xi}(\gamma, \mathcal{L}\varphi_t), 0 \le t \le 1\}$, then tests based on any continuous functionals of $\hat{\xi}_n(\gamma, \mathcal{L}\varphi_t)$ will be ADF for testing $H_0$. Some general sufficient conditions for the weak convergence of $\{\hat{\xi}_n(\gamma, \varphi_t), 0 \le t \le 1\}$ can be drawn from Proposition 6.2. Others can be inferred from, for example, van der Vaart and Wellner (1996). In particular, these claims hold for the family $\{\varphi_t, 0 \le t \le 1\}$ given at (3.10).

It is also important to note that the transformation $\mathcal{L}$ is free from $\gamma$ and, hence, the statement concerning the asymptotic distribution of $\{\hat{\xi}_n(\gamma, \mathcal{L}\varphi_t), 0 \le t \le 1\}$ is valid for any $\gamma \in L_2(\mathbb{R}^p, H)$.

Another consequence of Proposition 4.1 is worth formulating separately.

PROPOSITION 4.2. *Let* $\tilde{\theta}$ *be any estimator which satisfies* (2.3) *[and does not necessarily have a linear representation* (2.11)*] and let* $\psi_f$ *be a function of bounded variation. If, additionally,* (1.2) *and* (2.4) *hold, then for every* $\alpha = \gamma \cdot \varphi$ *with* $\gamma \in L_2(\mathbb{R}^p, H)$ *and* $\varphi \in \Phi$, *under* $H_0$,

$$\tilde{\xi}_n(\gamma, \mathcal{L}\varphi) \xrightarrow{d} w(\alpha),$$

*while under alternatives* (1.3),

$$\tilde{\xi}_n(\gamma, \mathcal{L}\varphi) \xrightarrow{d} w(\alpha) + \langle \mathcal{L}\alpha, a \rangle.$$

This proposition shows that although we used asymptotically linear representations (2.11) and (2.13) of $\hat{\theta}$ and $\tilde{\theta}$ to develop the previous theory, for the asymptotic behavior of the transformed processes $\hat{\xi}_n(\gamma, \mathcal{L}\varphi)$ and $\tilde{\xi}_n(\gamma, \mathcal{L}\varphi)$ the behavior of $\hat{\theta}$ and $\tilde{\theta}$ plays only a minor role.

It is instructive to consider informally a probabilistic connection between the processes $\hat{\xi}(\gamma, \varphi_t)$ and $\hat{\xi}(\gamma, \mathcal{L}\varphi_t)$. Let us associate with $\{\hat{\xi}(\gamma, \varphi_t), 0 \le t \le 1\}$ its natural filtration $\{\mathcal{F}_t, 0 \le t \le 1\}$, where each $\sigma$-field is

$$\widehat{\mathcal{F}}_t = \sigma\{\hat{\xi}(\gamma, \varphi_s), s \le t\}, \qquad 0 \le t \le 1,$$

and consider the filtered process $\{\hat{\xi}(\gamma, \varphi_t), \widehat{\mathcal{F}}_t, 0 \le t \le 1\}$. This is in $t$ a Gaussian semimartingale and it can be shown that the process $\{\hat{\xi}(\gamma, \mathcal{L}\varphi_t), \widehat{\mathcal{F}}_t, 0 \le t \le 1\}$ is actually its martingale part. In other words, if $\mathcal{V}$ denotes the Volterra operator defined by the integral on the right-hand side of (4.1), then the identity

$$(4.2) \qquad \hat{\xi}(\gamma, \varphi_t) = \hat{\xi}(\gamma, \mathcal{V}\varphi_t) + \hat{\xi}(\gamma, \mathcal{L}\varphi_t), \qquad 0 \le t \le 1,$$

is simply the Doob–Meyer decomposition of the process $\{\hat{\xi}(\gamma, \varphi_t), \widehat{\mathcal{F}}_t, 0 \le t \le 1\}$.

Details of this decomposition can be found in Khmaladze (1993), where the general construction of this form for a function-parametric process was introduced and studied. The notion of Doob–Meyer decomposition for a semimartingale can be found, for example, in Liptser and Shiryayev (1977).

REMARK 4.1. Since $\mathcal{L}\varphi$ is orthogonal to 1 and to $\psi_f$, the equality (4.2) can be rewritten in terms of the process $b$:

$$(4.3) \qquad b(\gamma, \varphi_t) = b(\gamma, \mathcal{V}\varphi_t) + b(\gamma, \mathcal{L}\varphi_t).$$

To some extent this is an unusual equation because both processes $b(\gamma, \varphi_t)$ and $b(\gamma, \mathcal{L}\varphi_t)$, taken separately, are Brownian motions. However, the nature of (4.3) can be more clearly understood as follows: let $\{\mathcal{F}_t^b, 0 \le t \le 1\}$ be the natural filtration of the process $b(\gamma, \varphi_t)$ in $t$ and let us enrich it with the $\sigma$-field $\sigma\{b(\gamma, h)\}$. Then the process $\{b(\gamma, \varphi_t), \mathcal{F}_t^b \vee \sigma\{b(\gamma, h)\}, 0 \le t \le 1\}$ is a Gaussian semimartingale (and not a martingale) and (4.3) is its Doob–Meyer decomposition. See, for example, Liptser and Shiryayev (1989) for more details on this.

REMARK 4.2. Another consequence of the orthogonality of $\mathcal{L}\varphi_t$ to 1 and to $\psi_f$ is this: although the process $\xi(\gamma, \varphi_t)$ with $\varphi_t$ chosen according to (3.10) with a nonconstant $\varphi$ is not a Brownian bridge (because in this case $\|\varphi_t^1\|^2 < \|\varphi_t\|^2 = t$) and hence even the process $\xi_n(\gamma, \varphi_t)$ with known value of parameter and statistics based on it may have an inconvenient limiting distribution, the transformed process $\xi(\gamma, \mathcal{L}\varphi_t)$ is the *standard Brownian motion* for any such choice of $\varphi_t$.

We shall now describe an analog of the above transformation suitable for testing the hypothesis $H_\sigma : G(y) = F(y/\sigma) \ \forall y \in \mathbb{R}$ and for some $\sigma > 0$. Let $\hat{\sigma}$ be an estimate of $\sigma$ based on $\{(X_i, Y_i), 1 \le i \le n\}$ satisfying

$$(4.4) \qquad \qquad \|n^{1/2}(\hat{\sigma} - \sigma)\| = O_p(1).$$

The analog of the processes $\hat{\xi}_n$ here is

$$\hat{\xi}_{n\sigma}(\gamma, \varphi) := n^{-1/2} \sum_{i=1}^{n} \gamma(X_i) \Big[ \varphi\Big( \frac{Y_i - \mu(X_i, \hat{\theta})}{\hat{\sigma}} \Big) - \int \varphi \, dF \Big].$$

To transform its weak limit $\hat{\xi}_\sigma$ under $H_\sigma$, again define an extended score function of $F((y - \mu)/\sigma)$ with respect to both parameters $\mu$ and $\sigma$, which is $h_\sigma(y) = (1, \psi_{f\mu}(y/\sigma), \psi_{f\sigma}(y/\sigma))^T$, where obviously

$$\psi_{f\mu}\Big( \frac{y}{\sigma} \Big) = \frac{1}{\sigma} \psi_f\Big( \frac{y}{\sigma} \Big),$$

$$\psi_{f\sigma}\Big( \frac{y}{\sigma} \Big) = \frac{1}{\sigma}\Big[ 1 + \frac{y}{\sigma} \psi_f\Big( \frac{y}{\sigma} \Big) \Big].$$

With notation

$$q(t) = \frac{1}{\sigma} f\Big( \frac{y}{\sigma} \Big),$$

$$q_\sigma(t) = \frac{y}{\sigma^2} f\Big( \frac{y}{\sigma} \Big), \qquad t = F\Big( \frac{y}{\sigma} \Big),$$

the analog of the $\Gamma_t$ matrix is

$$\Gamma_{\sigma,t} := \begin{pmatrix} 1 - t & -q(t) & -q_\sigma(t) \\ -q(t) & \int_t^1 \dot{q}^2(s) \, ds & \int_t^1 \dot{q}(s) \dot{q}_\sigma(s) \, ds \\ -q_\sigma(t) & \int_t^1 \dot{q}(s) \dot{q}_\sigma(s) \, ds & \int_t^1 \dot{q}_\sigma^2(s) \, ds \end{pmatrix}.$$

Again, assume that $\Gamma_{\sigma,t}^{-1}$ exists for all $0 \le t < 1$. Then, as above, let

$$(4.5) \quad \mathcal{L}_\sigma \varphi(y) := \varphi(y) - \int_{-\infty}^{y} \varphi(z) h_\sigma^T(z) \Gamma_{\sigma,F(z)}^{-1} \, dF(z) \, h_\sigma(y), \qquad y \in \mathbb{R}.$$

One can show that $\mathcal{L}_\sigma$ is a norm preserving transformation from $L_2(\mathbb{R}, F)$ to the subspace $\mathcal{H}_\sigma = \{\varphi \in L_2(\mathbb{R}, F) : \langle \varphi, h_\sigma \rangle = 0\}$ and hence $\hat{\xi}(\gamma, \mathcal{L}_\sigma \varphi)$ is a Brownian motion on $\mathbb{L}$.

PROOF OF PROPOSITION 4.1.  Though we could refer to the proof of Proposition 6.1, for presentational purposes it seems more convenient to give it here separately. Let, within this proof only, $\psi(t) := \varphi(F^{-1}(t))$ and $g(t) = h(F^{-1}(t))$ for $0 \leq t \leq 1$. Then

$$\int \mathcal{L}\varphi(y)h(y)^T \, dF(y)$$

$$= \int_0^1 \psi(t)g(t)^T \, dt - \int_0^1 \int_0^t \psi(s)g(s)^T \Gamma_s^{-1} \, ds \, g(t)g^T(t) \, dt$$

$$= \int_0^1 \psi(t)g(t)^T \, dt - \int_0^1 \psi(s)g(s)^T \Gamma_s^{-1} \int_s^1 g(t)g^T(t) \, dt \, ds$$

$$= \int_0^1 \psi(t)g(t)^T \, dt - \int_0^1 \psi(s)g(s)^T \, ds$$

$$= 0.$$

For the technical justification of the interchange of integration in the second equation above see the proof of Proposition 6.1 below or Khmaladze (1993). Similarly, we also have

$$\int [\mathcal{L}\varphi]^2 \, dF = \int_0^1 \psi^2(s) \, ds - 2 \int_0^1 \psi(s)g^T(s)\Gamma_s^{-1} \int_s^1 \psi(t)g(t) \, ds \, dt$$

$$+ \int_0^1 \int_0^1 \psi(s)g^T(s)\Gamma_s^{-1}\Gamma_{s \vee t}\Gamma_u h(t)\psi(t) \, ds \, dt$$

$$= \int \varphi^2 \, dF. \qquad \square$$

PROOF OF PROPOSITION 4.2. If $\psi_f$ is a function of bounded variation and $\varphi \in \Phi$, then $\mathcal{L}\varphi \in \Phi$ and therefore we can use (2.8). Together with the orthogonality property $\mathcal{L}\varphi \perp \psi_f$, which implies that $\langle \mathcal{L}\alpha, m_\theta \rangle = 0$, we obtain that

$$\hat{\xi}_n(\mathcal{L}\alpha) = \xi_n(\mathcal{L}\alpha) + o_p(1)$$

and the rest follows from Proposition 4.1, the CLT for $\xi_n(\mathcal{L}\alpha)$ and a standard contiguity argument. $\square$

## 5. Some explicit formulas and remarks.

5.1. *Transformation of the processes $\widehat{W}_1$ and $\hat{\xi}_n(1, \varphi_t)$.* In this section we shall apply the above transformation to residual empirical processes and give computational formulae of the transformed processes for testing $H_0$ and $H_\sigma$.

Recall from the previous sections that, for $0 \leq t \leq 1$,

$$(5.1) \qquad \widehat{U}_1(t) := \widehat{W}_1(F^{-1}(t)) = n^{-1/2} \sum_{i=1}^n [\mathbb{I}\{\varepsilon_i(\hat{\theta}) \leq F^{-1}(t)\} - t]$$

and

$$(5.2) \quad \hat{\xi}_n(1, \varphi_t) = n^{-1/2} \sum_{i=1}^n \left[ \varphi(\varepsilon_i(\hat{\theta})) \mathbb{I}\{\varepsilon_i(\hat{\theta}) \le F^{-1}(t)\} - \int_{y \le F^{-1}(t)} \varphi(y) \, dF(y) \right],$$

where in (5.2) $\varphi_t(y) = \varphi(y)\mathbb{I}(y \le F^{-1}(t))$. Note that $\widehat{U}_1(t)$ also corresponds to the $\hat{\xi}_n(1, \varphi_t)$, with $\varphi_t(y) = \mathbb{I}\{y \le F^{-1}(t)\}$. As another practically useful consequence of orthogonality of $\mathcal{L}\varphi$ to 1, we have the following equality:

$$\hat{\xi}_n(\gamma, \mathcal{L}\varphi) = n^{-1/2} \sum_{i=1}^n \gamma(X_i) \mathcal{L}\varphi(\varepsilon_i(\hat{\theta})) + o_p(1).$$

It means that we only need to construct transformations of random summands in (5.1) and (5.2). Introduce vector-functions

$$G(z) = \int_{y \le z} \Gamma_{F(y)}^{-1} h(y) \, dF(y),$$

$$J(z) = \int_{y \le z} \varphi(y) \Gamma_{F(y)}^{-1} h(y) \, dF(y), \qquad z \in \mathbb{R}.$$

Then the transformation $\mathcal{L}$ of (4.1) applied to $\widehat{U}_1$ of (5.1) gives

$$\widehat{w}_{n1}(t) = n^{-1/2} \sum_{Pi=1}^n \left[ \mathbb{I}\{\varepsilon_i(\hat{\theta}) \le z\} - [1, \psi_f(\varepsilon_i(\hat{\theta}))] G(z \wedge \varepsilon_i(\hat{\theta})) \right],$$
$$(5.3) \qquad\qquad\qquad\qquad\qquad\qquad\qquad\qquad t = F(z),$$

while the transformation of (5.2) is

$$\widehat{w}_{n2}(t) = n^{-1/2} \sum_{i=1}^n \left[ \varphi(\varepsilon_i(\hat{\theta})) \mathbb{I}\{\varepsilon_i(\hat{\theta}) \le z\} - [1, \psi_f(\varepsilon_i(\hat{\theta}))] J(z \wedge \varepsilon_i(\hat{\theta})) \right],$$
$$(5.4) \qquad\qquad\qquad\qquad\qquad\qquad\qquad\qquad t = F(z).$$

Similarly, to describe ADF tests for $H_\sigma$ based on the analog of $\widehat{w}_1$, let now $\hat{r}_i \equiv \varepsilon_i(\hat{\theta})/\hat{\sigma}$, and let us consider the processes

$$n^{-1/2} \sum_{i=1}^n [\mathbb{I}(\hat{r}_i \le z) - F(z)], \qquad n^{-1/2} \sum_{i=1}^n \left[ \varphi(\hat{r}_i) \mathbb{I}\{\hat{r}_i \le z\} - \int_{y \le z} \varphi(y) \, dF(y) \right],$$
$$t = F(z), \ z \in \mathbb{R}.$$

Then arguing as above, we are led to the following respective computational formulae:

$$w_{n1}(t) = n^{-1/2} \sum_{i=1}^n [\mathbb{I}\{\hat{r}_i \le z\} - h_{\hat{\sigma}}^T(\varepsilon_i(\hat{\theta})) G_{\hat{\sigma}}(z \wedge \varepsilon_i(\hat{\theta}))],$$

$$w_{n2}(t) = n^{-1/2} \sum_{i=1}^n [\varphi(\hat{r}_i) \mathbb{I}\{\hat{r}_i \le z\} - h_{\hat{\sigma}}^T(\varepsilon_i(\hat{\theta})) J_{\hat{\sigma}}(z \wedge \varepsilon_i(\hat{\theta}))],$$
$$t = F(z), \ z \in \mathbb{R},$$

where $h_\sigma(y)$ is as in the previous section, while $G_\sigma$ and $J_\sigma$ are defined as with $h$ replaced by $h_\sigma$ and $\Gamma$ replaced by $\Gamma_\sigma$.

These formulae may be used in the computation of any test statistic based on continuous functionals of $w_{n1}$, $w_{n2}$. From the theory developed above, if these functionals are invariant under the usual time transformation $t = F(y)$, they will be ADF!

5.2. *Nonrandom design.* We now state some analogous facts for the case of a nonrandom design where now the design vectors are denoted by $x_{ni}$. An analog of the condition (2.4) here is as follows: There exist a $q$-vector $\dot\mu$ on $\mathbb{R}^p \times \Theta$ and a $q \times q$ positive definite symmetric matrix $\Sigma$ such that

$$\Sigma_n := n^{-1} \sum_{i=1}^n \dot\mu(x_{ni}, \theta)\dot\mu^T(x_{ni}, \theta) \to \Sigma,$$

(5.5) $$\max_{1 \le i \le n} n^{-1/2}\|\dot\mu(x_{ni}, \theta)\| = o(1),$$

$$\sup_{1 \le i \le n, n^{1/2}\|\vartheta - \theta\| \le k} n^{1/2}|\mu(x_{ni}, \vartheta) - \mu(x_{ni}, \theta) - (\vartheta - \theta)^T\dot\mu(x_{ni}, \theta)| = o(1).$$

Under these conditions on the regression function and the rest of the conditions as before, the analogs of the above results with $\mu(X_i, \cdot)$ replaced by $\mu(x_{ni}, \cdot)$ remain valid in the present case. Using the results from Koul (1996), it is possible to obtain the analog of the expansions (2.14) and (2.15) under more general conditions on the function $\mu$ than given in (5.5), but we refrain from doing this for the sake of not obscuring main ideas and for the sake of brevity.

A similar remark applies to the linear regression model. In particular, in the case of nonrandom and general designs, but having the $n \times p$ design matrix $\mathbf{X}$ of rank $p$, just replace $n^{-1/2}X_i$ in the above formulas by $(\mathbf{X}'\mathbf{X})^{-1/2}x_{ni}$, $1 \le i \le n$, everywhere. Then tests based on the analogues of $\hat{w}_{n1}$ and $\hat{w}_{n2}$ are ADF for $H_0$, provided $\max_{1 \le i \le n} n^{1/2}\|(\mathbf{X}'\mathbf{X})^{-1/2}x_{ni}\| = O(1)$.

5.3. *Autoregressive time series.* Because of the close connection between regression and autoregressive models, analogues of the above ADF tests pertaining to the error distribution are easy to see in this case. Accordingly, suppose $Y_i$, $i \in \mathbb{Z} := \{0, \pm 1, \pm 2, \dots\}$, is now an observable stationary and ergodic time series. Let $\mu$ be as before satisfying (1.1) with $X_i := (Y_{i-1}, \dots, Y_{i-p})^T$, where $p \ge 1$ is a known integer. Then the above tests with this $X_i$ will be again ADF for testing $H_0$. A rigorous proof of this claim is similar to that appearing above, with the proviso that one uses the ergodic theorem in place of the law of the large numbers, and the CLT for martingale differences in place of the Lindeberg–Feller CLT. Note that now $H$ is the d.f. of the random vector $X_0$.

In the case of a stationary and ergodic *linear* AR($p$) model, that is, when $\mu(x, \vartheta) = x'\vartheta$, if the null error d.f. $F$ has mean zero and finite variance, then $EX_0 = 0$, that is, (3.5) is automatically satisfied, and hence tests based on the analog of $\widehat{U}_1$ of (5.1) will be a priori ADF for $H_0$. This was first proved in Boldin (1982), assuming $F$ has bounded second derivative, and in Koul (1991) when $F$ has only a uniformly continuous density. Thus, in linear autoregressive models the above transformation is useful only when there is a nonzero mean present in these models.

**6. Fitting a regression model.** In this section we shall develop some tests based on innovation processes that will be asymptotically distribution free for fitting a parametric model to the regression function $m(x) := E(Y|X = x)$. Actually we consider a somewhat more general problem where we fit a parametric model to a general regression function defined as follows.

For a real-valued measurable function $\varphi$ on $\mathbb{R}$, let $\mathcal{F}_\varphi$ denote a class of distribution functions $F$ on $\mathbb{R}$ such that $\varphi \in L_2(\mathbb{R}, F)$ and $\int |\varphi(y+t)| F(dy) < \infty$ for all $|t| \le k < \infty$. Let $m_\varphi(x)$ be defined by the relation

(6.1) $$E[\varphi(Y - m_\varphi(x))|X = x] = 0.$$

Note that if $\varphi(y) = y$, then $m_\varphi(x) = m(x)$, while if $\varphi(y) \equiv \mathbb{I}\{y > 0\} - (1 - \alpha)$, for an $0 < \alpha < 1$, then $m_\varphi(x)$ is the $\alpha$th quantile of the conditional distribution of $Y$, given $X = x$. The choice of $\varphi$ is up to the practitioner. The d.f. $F$ of the error $Y - m_\varphi(X)$ will be assumed to be an unknown member of $\mathcal{F}_\varphi$ for a given $\varphi$.

The problem of testing $\widetilde{H}_0$ is now extended to testing the hypothesis that $H_\varphi : m_\varphi(x) = \mu(x, \theta)$ for some $\theta \in \Theta$ against the alternatives described in (1.6). Consider again the function-parametric regression process

$$\xi_n(\gamma, \varphi; \vartheta) := n^{-1/2} \sum_{i=1}^{n} \gamma(X_i)\varphi(Y_i - \mu(X_i, \vartheta)).$$

Note that because of (6.1), under $H_\varphi$ $E\xi_n(\gamma, \varphi; \theta) = 0$.

Let $\tilde{\theta}$ be an M-estimator of $\theta$ satisfying (2.9) corresponding to $\eta_\vartheta \equiv \dot{\mu}_\vartheta$. Suppose, additionally, $F \in \mathcal{F}_\varphi$ is such that the function $t \mapsto \int \varphi(y+t)F(dy)$, $t \in \mathbb{R}$, is strictly monotonic and differentiable in a neighborhood of 0. Now, if we consider problems where $\varphi(y)$ is differentiable, such as $\varphi(y) = y$, which is a most interesting case, then we need to assume regularity condition (2.4) on the regression function $\mu(\cdot, \vartheta)$. While in the case of a nondifferentiable $\varphi$, as in, for example, $\varphi(y) = \mathbb{I}\{y > 0\} - (1 - \alpha)$, we need to assume as well that $F$, although unknown, satisfies also (1.2). In both cases, under (2.4) and (2.10) $\tilde{\theta}$ satisfies (2.11) and we obtain

$$\tilde{\xi}_n(\gamma, \varphi) = \xi_n(\gamma, \varphi) - \langle \gamma, \dot{\mu}_\theta^T \rangle C_\theta^{-1} \xi_n(\dot{\mu}_\theta, \varphi) + o_p(1)$$
$$= \xi_n(\gamma_\perp, \varphi) + o_p(1),$$

where

$$\gamma_\perp(x) = \gamma(x) - \langle \gamma, \dot{\mu}_\theta^T \rangle C_\theta^{-1} \dot{\mu}_\theta(x), \qquad x \in \mathbb{R}^p,$$

is the part of $\gamma$ orthogonal to $\dot{\mu}_\theta$ and no transformation of $\varphi$ is involved.

We emphasize that it is only for motivational purposes we are confining attention here to M-estimators. As we shall see later, any $n^{1/2}$-consistent estimator may be used to construct ADF tests for $H_\varphi$.

Now one can show that under $H_\varphi$, for each $\gamma, \varphi$ of the given type,

(6.2) $$\tilde{\xi}_n(\gamma, \varphi) \xrightarrow{d} b(\gamma_\perp, \varphi),$$

while under any sequence of alternatives (1.6),

$$\tilde{\xi}_n(\gamma, \varphi) \xrightarrow{d} b(\gamma_\perp, \varphi) + \lambda \langle \gamma_\perp, \ell_\theta \rangle,$$

where $\lambda$ is either $\langle \varphi', 1 \rangle$ or $-\langle \varphi, \psi_f \rangle$ depending on whether we assume (2.4) and (2.5) or (1.2).

As this last result shows, the asymptotic shift of the regression process $\tilde{\xi}_n(\gamma, \varphi)$ under the alternatives (1.6) is the linear functional of $\ell_\theta$ defined by the function $\gamma_\perp$. Therefore, to be able to detect all alternatives of the assumed type, we need to have a substantial supply of $\gamma_\perp$, that is, we need to consider $\tilde{\xi}_n(\gamma, \varphi)$ as a process in $\gamma$, and there is no need to vary $\varphi$ just in the same way as we had to vary $\varphi$ when testing our previous hypothesis $H_0$ and keep $\gamma$ fixed. We do not try to choose in any sense "optimal" $\varphi$ because the result will depend on $F$, while we prefer to work under the assumption that we do not know this d.f. Thus we can and will assume that $\varphi$ in the rest of this section is fixed.

From (6.2) we note that the limiting process as a function in $\gamma$ is again a projection of Brownian motion, but as a function in $\gamma_\perp$, it is just a Brownian motion.

Now we may have a convenient and customary way to parameterise $b(\gamma, \varphi)$ in $\gamma \in L_2(\mathbb{R}^p, H)$ to obtain processes with a standard and convenient distribution, and if we had similar ways to do this in subspaces of $L_2(\mathbb{R}^p, H)$, we could have the same convenient limiting processes in our problem. This, however, is not a straightforward task, as we have said earlier, especially because these subspaces, being orthogonal to $\dot{\mu}_\theta$, change from one regression function to another, and may even well change for the same regression function along the changes in the parameter $\theta$.

Nevertheless, we will see below that given a "convenient" indexing class $\mathcal{G}_0 \subset L_2(\mathbb{R}^p, H)$, in the sense that $\{b(\gamma, \varphi), \gamma \in \mathcal{G}_0\}$ forms a "convenient" asymptotic process—say, we can find the distribution of statistics based on $\{b(\gamma, \varphi), \gamma \in \mathcal{G}_0\}$ easily, and so on—we can map it isometrically into the subspace of functions orthogonal to $\dot{\mu}_\theta$. Thus, we obtain the process $\{b(\gamma, \varphi), \gamma \in \mathcal{G}_0'\}$, where $\mathcal{G}_0'$ is the image of this isometry, which on the one hand has exactly the same distribution and therefore carries the same "convenience" as the process $\{b(\gamma, \varphi), \gamma \in \mathcal{G}_0\}$, and on the other hand, is the limiting process for $\tilde{\xi}_n(\gamma, \varphi)$ if we index it by $\gamma \in \mathcal{G}'_0$.

To achieve this goal, first introduce the so called scanning family of measurable subsets $\mathcal{A} := \{A_z : z \in \mathbb{R}\}$ of $\mathbb{R}^p$ such that $A_z \subseteq A_{z'}$, for all $z \leq z'$, $H(A_{-\infty}) = 0$, $H(A_\infty) = 1$, and $H(A_z)$ is a strictly increasing absolutely continuous function of $z \in \mathbb{R}$.

To give examples, let $X^j$ denote the $j$th coordinate of the $p$-dimensional design variable $X$, $j = 1, \ldots, p$. Suppose that the marginal distribution of $X^1$ is absolutely continuous. Then we can take the family $A_z = \{x \in \mathbb{R}^p : x^1 \leq z\}$ as a scanning family. Or, if the sum $X^1 + \cdots + X^p$ is absolutely continuous, then one can take the family of half spaces $A_z = \{x \in \mathbb{R}^p : x^1 + x^2 + \cdots + x^p \leq z\}$.

Now let $B^c$ denote the complement of the set $B$,

$$z(x) := \inf\{z : A_z \ni x\},$$

$$\mathcal{C}_{\vartheta,z} := \int_{A_z^c} \dot{\mu}_\vartheta(y) \dot{\mu}_\vartheta^T(y) \, dH(y), \qquad\qquad z \in \mathbb{R},$$

$$\mathcal{T}_\vartheta \gamma(x) = \int_{A_{z(x)}} \gamma(y) \dot{\mu}_\vartheta^T(y) \mathcal{C}_{z(y)}^{-1} \, dH(y) \, \dot{\mu}_\vartheta(x), \qquad x \in \mathbb{R}^p, \ \vartheta \in \Theta.$$

We shall often write $\mathcal{C}_z$, $\mathcal{T}$ for $\mathcal{C}_{\theta,z}$, $\mathcal{T}_\theta$, respectively. Now, define the operator

$$\mathcal{K}\gamma(x) := \gamma(x) - \mathcal{T}\gamma(x), \qquad x \in \mathbb{R}^p.$$

PROPOSITION 6.1. *Let $\mathcal{G} := \{\gamma \in L_2(\mathbb{R}^p, H) : \langle \gamma, \dot{\mu}_\theta \rangle = 0\}$. Assume $\mathcal{C}_z$ is nonsingular for all $-\infty < z < \infty$. Then the transformation $\mathcal{K}$ is a norm preserving transformation from $L_2(\mathbb{R}^p, H)$ to $\mathcal{G}$:*

$$\mathcal{K}\gamma \perp \dot{\mu}_\theta, \qquad \|\mathcal{K}\gamma\| = \|\gamma\|.$$

*Consequently, for any fixed $\varphi$, the process $w(\gamma, \varphi) = \tilde{\xi}(\mathcal{K}\gamma, \varphi)$ is (function parametric) Brownian motion in $\gamma$.*

Similarly to Proposition 4.2, the following corollary shows that much less is required from an estimator $\tilde{\theta}$ than its asymptotic linearity. The random vector $Z$ below can be thought of as the limit in distribution of $\sqrt{n}(\tilde{\theta} - \theta)$.

COROLLARY 6.1. *Let $\tilde{\xi}$ be any process of the form*

$$\tilde{\xi}(\gamma, \varphi) = b(\gamma, \varphi) - \langle \gamma, \dot{\mu}_\theta^T \rangle Z,$$

*where $Z$ is a random vector (not necessarily Gaussian) in $\mathbb{R}^q$. Then for any fixed $\varphi \in L_2(\mathbb{R}, F)$, the process*

$$w(\gamma, \varphi) = \tilde{\xi}(\mathcal{K}\gamma, \varphi)$$

*is Brownian motion in $\gamma \in L_2(\mathbb{R}^p, H)$.*

Now we shall, as an example, focus on the case $\gamma = \mathbb{I}_B$, for $B$ a Borel set in $\mathbb{R}^p$. Then

$$\mathcal{K}\mathbb{I}_B(x) = \mathbb{I}_B(x) - \int_{A_{z(x)}} \mathbb{I}_B(y) \dot{\mu}_\theta^T(y) \mathcal{C}_{z(y)}^{-1} \, dH(y) \, \dot{\mu}_\theta(x).$$

In view of the above discussion, our transformation is the process

$$w_n(B) := \tilde{\xi}_n(\mathcal{K}\mathbb{I}_B, \varphi)$$

(6.3)
$$= n^{-1/2} \sum_{i=1}^n \left[ \mathbb{I}_B(X_i) - \int_{A_{z(X_i)}} \mathbb{I}_B(y) \dot{\mu}_\theta^T(y) \mathcal{C}_{z(y)}^{-1} \, dH(y) \, \dot{\mu}_\theta(X_i) \right]$$
$$\times \varphi(Y_i - \mu(X_i, \tilde{\theta})).$$

We do not consider in this paper the problem of weak convergence of transformed processes $\{\tilde{\xi}_n(\gamma, \mathcal{L}\varphi), \varphi \in \Phi_0\}$ or $\{\tilde{\xi}_n(\mathcal{K}\gamma, \varphi), \gamma \in \mathcal{G}_0\}$ to corresponding Brownian motions for appropriate indexing classes $\Phi_0$ and $\mathcal{G}_0$ in full generality. Nevertheless we shall now state a sufficient condition under which the process (6.3) converges weakly to a set-parametric Brownian motion on the practically important class of sets—a subclass $\mathcal{B}_0$ of all right closed rectangles in $\mathbb{R}^p$, that is, $\mathcal{B}_0 \subset \{(-\infty, v], v \in \mathbb{R}^p\}$. Our assumption is the following:

(6.4)        There exists a $\tau > 0$ such that $B \subseteq A_{1-\tau}$ for all $B \in \mathcal{B}_0$.

This condition is not necessary, but simplifies the proof substantially. See Khmaladze (1993) for the version without this condition.

Let $\{w(B), B \in \mathcal{B}_0\}$ be set-parametric Brownian motion on $\mathcal{B}_0$ with covariance function

$$Ew(B)w(B') = cH(B \cap B'),$$

where, without loss of generality we can assume the constant $c$ to be 1; compare Remark 6.1.

The space in which we will consider weak convergence of $w_n$ will be $\ell^\infty(\mathcal{G}_0)$, where $\mathcal{G}_0 = \{\mathbb{I}_B(\cdot), B \in \mathcal{B}_0\}$ is equipped with the $L_2$-norm. [See, e.g., page 34 in van der Vaart and Wellner (1996).] Now, write $\hat{\varepsilon}_i, \tilde{\varepsilon}_i$ for $\varepsilon_i(\hat{\theta})$, $\varepsilon_i(\tilde{\theta})$, respectively. Also, let $\varepsilon$ denote a r.v. having the same distribution as $\varepsilon_1(\theta)$.

PROPOSITION 6.2.    *Suppose regularity conditions* (2.4) *and* (2.5) *are satisfied. Suppose also* $E\varphi^2(\varepsilon) = 1$ *and* $\tilde{\theta}$ *is any estimator such that* $\sqrt{n}(\tilde{\theta} - \theta) = O_p(1)$. *If* $\mathcal{B}_0$ *is such that* (6.4) *is satisfied then, under* $H_\varphi$

$$w_n \xrightarrow{d} w, \qquad \text{in } l^\infty(\mathcal{G}_0).$$

REMARK 6.1.  In the definition (6.3) of the process $w_n$ we assumed that $E\varphi^2(\varepsilon) = 1$ without loss of generality. Indeed, we can always replace $\varphi(\hat{\varepsilon}_i)$ by $\varphi(\hat{\varepsilon}_i)/\hat{\sigma}$ in $w_n$, where $\hat{\sigma}^2 = n^{-1}\sum_{i=1}^n \varphi^2(\hat{\varepsilon}_i)$ is an estimator of $\sigma^2 = E\varphi^2(\varepsilon)$. Then it is obvious that the processes which incorporate $\varphi(\hat{\varepsilon}_i)/\hat{\sigma}$ and $\varphi(\hat{\varepsilon}_i)/\sigma$, respectively, will converge to each other, uniformly in $B$, in probability.

Since the kernel of the transformation $\mathcal{T}$ depends on $\theta$, we will certainly need to replace it with an estimator. It seems the simplest to use the same estimator $\hat{\theta}$ as is used in $\tilde{\xi}_n$, although it is not necessary and in principle any consistent estimator can be used: small perturbation of $\theta$ in $\mathcal{T}_\theta$ will only slightly perturb the process $\tilde{\xi}_n(\mathcal{T}_\theta\gamma, \varphi)$. To prove this latter statement formally, we need to complement (2.4) by the following two mild assumptions. Let

$$d^2(\vartheta_1, \vartheta_2) := E\|\dot{\mu}_{\vartheta_1}(X) - \dot{\mu}_{\vartheta_2}(X)\|^2 E\varphi^2(\tilde{\varepsilon}), \qquad \tilde{\varepsilon} := \varepsilon(\bar{\theta}), \ \vartheta \in \Theta,$$

$$\rho(\delta) := \sup_{\|\vartheta_1 - \vartheta_2\| \le \delta} d(\vartheta_1, \vartheta_2), \qquad\qquad \delta > 0.$$

Suppose that $E\varphi^2(\varepsilon) = 1$, and that for some $\epsilon > 0$,

$$(6.5) \qquad \sup_{\|\vartheta_1 - \vartheta_2\| \le \epsilon} \left| \frac{1}{n} \sum_{i=1}^n \|\dot{\mu}_{\vartheta_1}(X_i) - \dot{\mu}_{\vartheta_2}(X_i)\|^2 - d^2(\vartheta_1, \vartheta_2) \right| = o_p(1),$$

$$\text{as } n \to \infty,$$

$$(6.6) \qquad\qquad \sum_{k=0}^\infty k\rho(\epsilon 2^{-k}) < \infty.$$

Define the estimated tranformed process:

$$\widetilde{w}_n(B) := \tilde{\xi}_n(\mathbb{I}_B, \varphi) - \tilde{\xi}_n(\mathcal{T}_{\hat{\theta}}\mathbb{I}_B, \varphi).$$

We have the following statement.

PROPOSITION 6.3.  *Let $\{\mathbb{I}_B, B \in \mathcal{B}_0\}$ be any collection of indicator functions such that $\mathcal{B}_0$ satisfies (6.4). Then under the assumptions (6.5) and (6.6),*

$$\sup_{B \in \mathcal{B}_0} |\widetilde{w}_n(B) - w_n(B)| = o_p(1).$$

To prove this last proposition we will use the following lemma, which is of independent interest. Let, for a $c > 0$,

$$D_n = \left\{ \sup_{\|\vartheta_1 - \vartheta_2\| \le \delta} \frac{1}{n} \sum_{i=1}^n \|\dot{\mu}_{\vartheta_1}(X_i) - \dot{\mu}_{\vartheta_2}(X_i)\|^2 \varphi^2(\tilde{\varepsilon}_i) \le (1+c)d^2(\delta) \right.$$

$$\left. \text{for all } 0 < \delta < c \right\}.$$

LEMMA 6.1. *Let* $\{\mathbb{I}_A, A \in \mathcal{A}'\}$ *be any collection of indicator functions. Then under the assumptions* (6.5) *and* (6.6),

$$P\left(\sup_{\|\vartheta - \theta\| \leq \epsilon} |\tilde{\xi}_n(\mathbb{I}_A \dot{\mu}_\vartheta, \varphi) - \xi_n(\mathbb{I}_A \dot{\mu}_\theta, \varphi)| > x|D_n\right)$$

$$\leq \exp\left\{-(x/2)C\sum_{k=0}^{\infty} k\rho(\epsilon\, 2^{-k})\right\},$$

$$E\left\{\sup_{\|\vartheta - \theta\| \leq \epsilon} |\tilde{\xi}_n(\mathbb{I}_A \dot{\mu}_\vartheta, \varphi) - \tilde{\xi}_n(\mathbb{I}_A \dot{\mu}_\theta, \varphi)|^2|D_n\right\}$$

$$\leq C\sum_{k=0}^{\infty} k\rho(\epsilon\, 2^{-k}) \to 0,$$

*as* $\epsilon \to 0$, *where* $C$ *is a positive universal constant.*

Now we prove all three propositions and the lemma.

PROOF OF PROPOSITION 6.1. Fix a $k < \infty$ and consider $\gamma_k := \gamma\mathbb{I}_{A_k}$. We shall first show that $\langle \mathcal{K}\gamma_k, \dot{\mu}_\theta^T \rangle = 0$. Note that $y \in A_{z(x)}$ is equivalent to $x \in A_{z(y)}^c$ for almost all $x, y$ with respect to the measure $H$. This fact, together with changing the order of integration, yields

$$\langle \mathcal{K}\gamma_k, \dot{\mu}_\theta^T \rangle$$

$$= \int_{\mathbb{R}^p} \gamma_k(x)\dot{\mu}_\theta^T(x)\,dH(x)$$

$$\quad - \int_{\mathbb{R}^p} \int_{A_{z(x)}} \gamma_k(y)\dot{\mu}_\theta^T(y)\mathcal{C}_{z(y)}^{-1}\,dH(y)\dot{\mu}_\theta(x)\dot{\mu}_\theta^T(x)\,dH(x)$$

$$= \langle \gamma_k, \dot{\mu}_\theta^T \rangle - \int_{\mathbb{R}^p} \gamma_k(y)\,\dot{\mu}_\theta^T(y)\mathcal{C}_{z(y)}^{-1}\,dH(y)\int_{A_{z(y)}^c} \dot{\mu}_\theta(x)\dot{\mu}_\theta^T(x)\,dH(x)$$

$$= \langle \gamma_k, \dot{\mu}_\theta^T \rangle - \langle \gamma_k, \dot{\mu}_\theta^T \rangle$$

$$= 0.$$

Now we shall show that

$$(6.7) \qquad\qquad \langle \mathcal{K}\gamma_k, \mathcal{K}\gamma_k \rangle = \langle \gamma_k, \gamma_k \rangle.$$

Using the notation

$$\rho_k^T(z) := \int_{A_z} \gamma_k(y)\dot{\mu}_\theta^T(y)\mathcal{C}_{z(y)}^{-1}\,dH(y), \qquad z \in \overline{\overline{\mathbb{R}}},$$

rewrite

$$\langle \mathcal{K}\gamma_k, \mathcal{K}\gamma_k \rangle$$
$$= \langle \gamma_k, \gamma_k \rangle - 2 \int_{\mathbb{R}^p} \rho_k^T(z(x)) \dot{\mu}_\theta(x) \gamma(x) \, dH(x)$$
$$+ \int_{\mathbb{R}^p} \rho_k^T(z(x)) \dot{\mu}_\theta(x) \dot{\mu}_\theta^T(x) \rho(z(x)) \, dH(x)$$
$$= \langle \gamma_k, \gamma_k \rangle - 2 \int_{z \leq z_0} \rho_k^T(z) \mathcal{C}_z \, d\rho_k(z) + \int_{-\infty}^{\infty} \rho_k^T(z) \, d\mathcal{C}_z \rho_k(z)$$
$$= \langle \gamma_k, \gamma_k \rangle - \rho_k^T(z) \mathcal{C}_z \rho_k(z)|_{-\infty}^{\infty}.$$

Because $\gamma_k = \gamma \mathbb{I}_{A_k}$, the function $\rho_k$ remains bounded as $z \to \infty$ and hence the substitution in the above equals zero, thereby proving (6.7).

Next, by definition $\gamma_k \to \gamma$ as $k \to \infty$. Let $k \to \infty$ in (6.7) to conclude that it remains true for a general $\gamma \in L_2(\mathbb{R}^p, H)$.   □

PROOF OF PROPOSITION 6.2.   Using the definition of the operator $\mathcal{K}$ one can write

$$\sup_{B \in \mathcal{B}_0} |w_n(B) - \xi_n(\mathcal{K}\mathbb{I}_B, \varphi)|$$
$$\leq \sup_{B \in \mathcal{B}_0} |\tilde{\xi}_n(\mathbb{I}_B, \varphi) - \xi_n(\mathbb{I}_B, \varphi) - E\varphi' E(\mathbb{I}_B \dot{\mu}_\theta^T) \, n^{1/2}(\tilde{\theta} - \theta)|$$
$$+ \sup_{B \in \mathcal{B}_0} |\tilde{\xi}_n(\mathcal{T}\mathbb{I}_B, \varphi) - \xi_n(\mathcal{T}\mathbb{I}_B, \varphi) + E\varphi' E(\mathbb{I}_B \dot{\mu}_\theta^T) n^{1/2}(\tilde{\theta} - \theta)|.$$

However, Proposition 2.1 implies that the first supremum on the right-hand side is $o_p(1)$. To deal with the second supremum, let us use the fact that $\mathbb{I}_{A_{z(x)}}(y) = \mathbb{I}_{A_{z(y)}^c}(x)$ a.e. and change the order of summation and integration:

$$\xi_n(\mathcal{T}\gamma, \varphi) = n^{-1/2} \sum_{i=1}^n \int_{A_t(X_i)} \gamma(y) \, \dot{\mu}_\theta^T(y) C_{z(y)}^{-1} \, dH(y) \dot{\mu}_\theta(X_i) \varphi(e_i)$$
$$= \int \gamma(y) \dot{\mu}_\theta^T(y) C_{z(y)}^{-1} \xi_n(\mathbb{I}_{A_t^c(y)} \dot{\mu}_\theta, \varphi) \, dH(y).$$

Similar equality is certainly true for $\tilde{\xi}_n$. Therefore, using Proposition 2.1 once again, we obtain

$$\sup_{B \in \mathcal{B}_0} |\tilde{\xi}_n(\mathcal{T}\mathbb{I}_B, \varphi) - \xi_n(\mathcal{T}\mathbb{I}_B, \varphi) - E\varphi' E(\mathbb{I}_B \dot{\mu}_\theta^T) n^{1/2}(\tilde{\theta} - \theta)| = o_p(1).$$   □

PROOF OF PROPOSITION 6.3.   First note that from the previous proof we have

$$\tilde{w}_n(B) - w_n(B) = \tilde{\xi}_n(\mathcal{T}_{\tilde{\theta}} \mathbb{I}_B, \varphi) - \tilde{\xi}_n(\mathcal{T}_\theta \mathbb{I}_B, \varphi).$$

Now let

$$\eta^T(y, \vartheta) = \dot{\mu}_\vartheta^T(y) \mathcal{C}_{\vartheta, z(y)}^{-1}$$

and

$$\tilde{\xi}_n(z, \vartheta) = \tilde{\xi}_n(\mathbb{I}_{A_{z(y)}^c} \dot{\mu}_\vartheta, \varphi).$$

Then we can rewrite

$$\tilde{\xi}_n(\mathcal{T}_\vartheta \mathbb{I}_B, \varphi) = \int \mathbb{I}_B(y) \dot{\mu}_\vartheta^T(y) \mathcal{C}_{\vartheta, z(y)}^{-1} \tilde{\xi}_n(\mathbb{I}_{A_{z(y)}^c} \dot{\mu}_\vartheta, \varphi) \, dH(y)$$

$$= \int \mathbb{I}_B(y) \eta^T(y, \vartheta) \tilde{\xi}_n(z(y), \vartheta) \, dH(y).$$

Since

$$|\tilde{\xi}_n(\mathcal{T}_{\hat{\theta}} \mathbb{I}_B, \varphi) - \tilde{\xi}_n(\mathcal{T}_\theta \mathbb{I}_B, \varphi)|$$

$$\leq \mathbb{I}_{\{\|\tilde{\theta} - \theta\| \leq \epsilon\}} \sup_{\|\vartheta - \theta\| \leq \epsilon} |\tilde{\xi}_n(\mathcal{T}_\vartheta \mathbb{I}_B, \varphi) - \tilde{\xi}_n(\mathcal{T}_\theta \mathbb{I}_B, \varphi)|$$

$$+ \mathbb{I}_{\{\|\tilde{\theta} - \theta\| > \epsilon\}} |\tilde{\xi}_n(\mathcal{T}_{\hat{\theta}} \mathbb{I}_B, \varphi) - \tilde{\xi}_n(\mathcal{T}_\theta \mathbb{I}_B, \varphi)|$$

and $\tilde{\theta}$ is consistent estimator, it is enough to prove that

$$\sup_{B \in \mathcal{B}_0, \|\vartheta - \theta\| \leq \epsilon} \int \mathbb{I}_B(y) |\eta^T(y, \vartheta) \tilde{\xi}_n(z(y), \vartheta) - \eta^T(y, \theta) \tilde{\xi}_n(z(y), \theta)| \, dH(y)$$

$$= o_p(1),$$

as $\epsilon \to 0$ and $n \to \infty$. Using the Cauchy–Schwarz inequality and the fact that $B \subset A_{1-\tau}$, we find that the left-hand side of the above equality is bounded above by

$$\sup_{\|\vartheta - \theta\| \leq \epsilon} \|\eta(\cdot, \vartheta) - \eta(\cdot, \theta)\|_H \|\tilde{\xi}_n(\cdot, \vartheta)\|_H$$

$$+ \|\eta(\cdot, \theta)\|_H \sup_{\|\vartheta - \theta\| \leq \epsilon} \|\tilde{\xi}_n(\cdot, \vartheta) - \tilde{\xi}_n(\cdot, \theta)\|_H,$$

where $\| \cdot \|_H$ is the $L_2$ norm with respect to $H$.

Since $C_z$ is nonsingular for $z < 1 - \tau$, we have $\|\eta(\cdot, \theta)\|_H < \infty$. Moreover, $\dot{\mu}_\vartheta$ being continuous in $\vartheta$ in mean square sense [condition (2.4)], it follows that for all sufficiently small $\epsilon$, $C_{\vartheta, z}$ is nonsingular for all $\|\vartheta - \theta\| \leq \epsilon$, $z < 1 - \tau$, and that $\sup_{\|\vartheta - \theta\| \leq \epsilon} \|\eta(\cdot, \vartheta) - \eta(\cdot, \theta)\|_H$ is small. What remains therefore to show is that $\sup_{B \in \mathcal{B}_0, \|\vartheta - \theta\| \leq \epsilon} \|\tilde{\xi}_n(\cdot, \vartheta)\|_H = o_p(1)$, and that $\sup_{B \in \mathcal{B}_0, \|\vartheta - \theta\| \leq \epsilon} \|\tilde{\xi}_n(\cdot, \vartheta) - \tilde{\xi}_n(\cdot, \theta)\| = o_p(1)$ as $n \to \infty$ and $\epsilon \to 0$. These properties are proved in Lemma 6.1. □

PROOF OF LEMMA 6.1. First note that a symmetrization lemma [see, e.g., van der Vaart and Wellner (1996), Section 2.3.2] can be used to imply that

$$\|\tilde{\xi}_n(z, \vartheta_1) - \tilde{\xi}_n(z, \vartheta_2)\| \le 2\|\tilde{\xi}_n^0(z, \vartheta_1) - \tilde{\xi}_n^0(z, \vartheta_2)\|,$$

where

$$\tilde{\xi}_n^0(z, \vartheta) = n^{-1/2} \sum_{i=1}^n e_i \mathbb{I}_{A_z^c} \dot{\mu}_\vartheta(X_i) \varphi(\tilde{\varepsilon}_i),$$

and $\{e_i\}_{i=1}^n$ are Rademacher random variables independent of $\{(X_i, Y_i)\}_{i=1}^n$. Averaging first over $\{e_i\}_{i=1}^n$, we obtain for all $t > 0$,

$$E[\exp\{t^{-1}\|\tilde{\xi}_n^0(z, \vartheta_1) - \tilde{\xi}_n^0(z, \vartheta_2)\|\}|D_n]$$

$$\le E\left[\exp\left\{2t^{-2}n^{-1}\sum_{i=1}^n \|\dot{\mu}_{\vartheta_1}(X_i) - \dot{\mu}_{\vartheta_2}(X_i)\|^2 \varphi^2(\tilde{\varepsilon}_i)\right\}\Big|D_n\right]$$

$$\le \exp\{2t^{-2}(1+c)\rho^2(\|\vartheta_1 - \vartheta_2\|)\}.$$

Following van der Vaart and Wellner (1996), Section 2.2, denote by $\|X\|_{\psi, D_n}$ the Orlicz norm of the random variable $X$ induced by the function $\psi(x) = e^x - 1$—this is the smallest constant $t$ such that $E[\exp(|X|/t) - 1|D_n] \le 1$. Then the previous inequality implies that

$$\|\tilde{\xi}_n^0(z, \vartheta_1) - \tilde{\xi}_n^0(z, \vartheta_2)\|_{\psi, D_n} \le C\rho(\|\vartheta_1 - \vartheta_2\|).$$

Since $e^{x/t} - 1 > (x/t)^2/2!$, it immediately follows that

$$E \sup_{\|\vartheta - \theta\| \le \epsilon} \|\tilde{\xi}_n(y, \vartheta) - \tilde{\xi}_n(y, \theta)\|^2$$

$$\le 2\left\|\sup_{\|\vartheta - \theta\| \le \epsilon} \|\tilde{\xi}_n(y, \vartheta) - \tilde{\xi}_n(y, \theta)\|\right\|_{\psi, D_n}.$$

We now show that the Orlicz norm on the right-hand side is small for small $\epsilon$. We will do this by slightly adjusting the chaining argument. Let $N(\delta)$ be the covering number [the cardinality of the minimal $\delta$-net $\mathcal{N}(\delta)$] of the unit ball in $\mathbb{R}^q$. Let each $\zeta_{k+1} \in \mathcal{N}(2^{-k-1})$ be linked to unique $\zeta_k \in \mathcal{N}(2^{-k})$ in such a way that $\|\zeta_{k+1} - \zeta_k\| \le 2^{-k}$. Then using the Fundamental Lemma 2.2 of van der Vaart and Wellner (1996), Section 2.2, one can write (with $\vartheta = \epsilon\zeta$)

$$\left\|\max_{\vartheta_k, \vartheta_{k+1}} \|\tilde{\xi}_n(y, \vartheta_k) - \tilde{\xi}_n(y, \vartheta_{k+1})\|\right\|_{\psi, D_n}$$

$$\le C \ln N(2^{-k}) \rho(\epsilon 2^{-k}).$$

Hence

$$\left\| \sup_{\|\vartheta - \theta\| \leq \epsilon} \|\tilde{\xi}_n^0(y, \vartheta) - \tilde{\xi}_n^0(y, \theta)\| \right\|_{\psi, D_n}$$

(6.8)
$$\leq \sum_{k=1}^{\infty} \left\| \max_{\vartheta_k, \vartheta_{k+1}} \|\tilde{\xi}_n(y, \vartheta_k) - \tilde{\xi}_n(y, \vartheta_{k+1})\| \right\|_{\psi, D_n}$$

$$\leq C \sum_{k=1}^{\infty} \ln N(2^{-k}) \rho(\epsilon 2^{-k}) \leq Cq \sum_{k=1}^{\infty} k \rho(\epsilon 2^{-k}),$$

where the last inequality follows from obvious estimation from above, $N(\delta) \leq C\delta^{-q}$. Since $\rho(\epsilon 2^{-k}) \to 0$ as $\epsilon \to 0$ and the series converges for some $\epsilon > 0$, it tends to 0 as $\epsilon \to 0$.

Finally, combine the symmetrization and Markov inequalities to obtain

$$P\left( \sup_{\|\vartheta - \theta\| \leq \epsilon} \|\tilde{\xi}_n(y, \vartheta) - \tilde{\xi}_n(y, \theta)\| > x | D_n \right)$$

$$\leq P\left( \sup_{\|\vartheta - \theta\| \leq \epsilon} \|\tilde{\xi}_n^0(y, \vartheta) - \tilde{\xi}_n^0(y, \theta)\| > \frac{x}{2} \Big| D_n \right)$$

$$\leq E\left[ \exp\left\{ t^{-1} \sup_{\|\vartheta - \theta\| \leq \epsilon} \|\tilde{\xi}_n^0(y, \vartheta) - \tilde{\xi}_n^0(y, \theta)\| \right\} \Big| D_n \right] \exp\left( -\frac{x}{4t} \right).$$

From the definition of the Orlicz norm $\| \sup_{\|\vartheta - \theta\| \leq \epsilon} \|\tilde{\xi}_n^0(y, \vartheta) - \tilde{\xi}_n^0(y, \theta)\| \|_{\psi, D_n}$ and the inequality (6.8), it follows that the expectation above does not exceed 2 for $t = q \sum_{k=1}^{\infty} k \rho(\epsilon 2^{-k})$. Hence the inequality of the lemma. $\quad \square$

We end this section by pointing out that the conditions (6.5) and (6.6) are trivially satisfied in the case $\mu(x, \vartheta) \equiv \vartheta' S(x)$, where $S(x)$ is a vector of functions of $x$ with finite second moment $E\|S(X)\|^2$.

## 7. Some simulations.

This section presents some simulations to see how well the finite sample level of significance is approximated by the asymptotic level for the supremum of the absolute values of the transformed processes defined at (5.3) and (6.3). It is noted that when fitting a standard normal distribution to the errors with a rapidly changing regression function, or when fitting a two-variable linear regression model with standard normal errors and using the least squares residuals, this approximation is very good even for the sample size 40, especially in the right tail.

The lack of an analytical form of the distribution of the supremum of the Brownian motion on $[0, 1]^2$ created an extra difficulty here. We had to first obtain simulated approximation to this distribution. This was done by simulating an appropriate two time parameter Poisson process of sample size $5K$, with $20K$ replications. Selected quantiles based on this simulation are presented in Tables 2 of Section 7.2. This should be of independent interest also.

7.1. $\sup_z |\widehat{w}_{n1}(z)|$ *of* (5.3). This section presents some selected empirical percentiles of the transformed statistic $D_n := \sup_z |\widehat{w}_{n1}(z)|$ of (5.3) for testing $H_0: F$ is the standard normal d.f. The regression function is taken to be $\mu(x, \vartheta) = e^{\vartheta x}$, with true $\theta = 0.25$, the regressors $X_i, i = 1, \ldots, n$, are chosen to be uniformly distributed on $[2, 4]$, and the errors $\varepsilon_i \equiv \varepsilon_i(\theta), i = 1, \ldots, n$, are standard Gaussian. In this case the $\Gamma$ function of Section 5.3 becomes

$$\Gamma_{F(y)}^{-1} = \frac{1}{[1 - F(y)][ya(y) + 1 - a^2(y)]} \begin{pmatrix} 1 + ya(y) & a(y) \\ a(y) & 1 \end{pmatrix},$$

where $a(y) = f(y)/(1 - F(y))$, with $f$ and $F$ denoting the standard normal density and d.f., respectively. Consequently, the vector-function $G$ of (5.3) is now equal to

$$\begin{aligned} G^T(z) &= \int_{-\infty}^{z} (1, -y) \Gamma_{F(y)}^{-1} f(y) \, dy \\ &= \int_{-\infty}^{z} \frac{1}{ya(y) + 1 - a^2(y)} (1, a(y) - y) a(y) \, dy \end{aligned}$$

and, eventually, the transformed process of (5.3) has the form

$$\begin{aligned} \widehat{w}_{n1}(t) = n^{-1/2} \sum_{i=1}^{n} \Big[ & \mathbb{I}\{\varepsilon_i(\hat{\theta}) \leq z\} \\ (7.1) \\ &- \int_{-\infty}^{z \wedge \varepsilon_i(\hat{\theta})} \frac{1 + \varepsilon_i(\hat{\theta})(a(y) - y)}{ya(y) + 1 - a^2(y)} a(y) \, dy \Big], \qquad t = F(z). \end{aligned}$$

Although the form of the regression function does not participate in the martingale transformation $\mathcal{L}$ it still may affect the finite sample behavior of the transformed process as far as it affects $\varepsilon_i(\hat{\theta}), i = 1, \ldots, n$, where $\hat{\theta}$ is the MLE under the null hypothesis. It was thus of interest to see whether the estimation of $\theta$ will not affect the values of $\varepsilon_i(\hat{\theta}), i = 1, \ldots, n$, too much and worsen the convergence of the transformed process to its limit. For this reason we chose a more or less rapidly changing regression function. On the other hand there was no point in choosing multidimensional regressors $X_i$ here, since the transformed process depends solely on $\varepsilon_i(\hat{\theta}), i = 1, \ldots, n$.

We simulated $\{(X_i, Y_i = e^{0.25X_i} + \varepsilon_i), 1 \leq i \leq n\}$ for sample sizes $n = 40$, 100 and for each sample calculated the value of the Kolmogorov–Smirnov statistic $D_n := \sup\{|\widehat{w}_{n1}(t)|; 0 \leq t \leq 1\}$, with $\widehat{w}_{n1}(t)$ as in (7.1). This was done $m = 10K$ times. In Table 1 $d_\alpha$ is the $100(1 - \alpha)\%$ percentile of the limiting distribution of $D_n$. The values are obtained by approximating the d.f. of the supremum of the Brownian motion over $[0,1]$ by $\mathcal{G}(z) := P(\sup_{0 \le t \le 1} |\xi_n(t) - nt|/\sqrt{n} \le z)$, with $n = 5K$, where $\xi_n(t)$, $t \in [0,1]$, is a Poisson process with intensity $n$. The d.f. $\mathcal{G}$ was calculated using the exact recurrence formulas and code given in Khmaladze and Shinjikashvili (2001). The values obtained are accurate to $5 \cdot 10^{-3}$.

TABLE 1
*Selected quantiles of $P(D_n > d_\alpha)$*

| $\alpha$ | 0.2 | 0.1 | 0.05 | 0.025 | 0.01 |
|---|---|---|---|---|---|
| $n \setminus d_\alpha$ | 1.64 | 1.96 | 2.24 | 2.50 | 2.81 |
| 40 | 0.168 | 0.084 | 0.046 | 0.029 | 0.019 |
| 100 | 0.178 | 0.093 | 0.052 | 0.029 | 0.014 |

Table 1 also gives the Monte Carlo estimates of $P(D_n > d_\alpha)$ for $n = 40$ and $n = 100$ based on $m = 10K$ replications. The resulting (simulated) distribution functions of $D_n$ along with $\mathcal{G}$ as solid line are shown in Figure 1. The quality of approximation appears to be quite close to what one has in the classical case of the empirical process and the limiting Brownian bridge especially in the upper tail, where we need it the most.

7.2. $\sup_B |w_n(B)|$ *of* (6.3). Here the regressors $X_i$, $i = 1, \dots, n$, are two-dimensional Gaussian random vectors with standard normal marginal distributions and correlation $r$. The regression function being fitted is chosen to be linear:

$$(7.2) \qquad \mu(x, \vartheta) = \vartheta_1 x_1 + \vartheta_2 x_2,$$

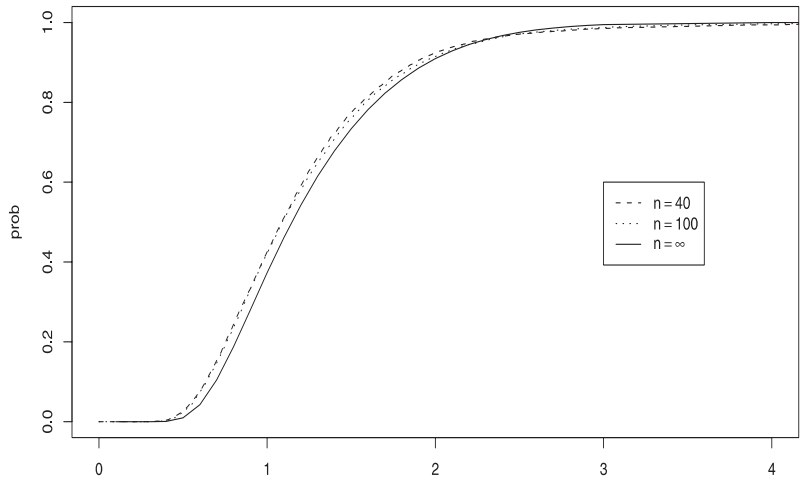

FIG. 1. *E.d.f. of $D_n$ for $n = 40, 100$, $m = 10K$ and $\mathcal{G}$.*

with the true parameter $\theta' = (1,1)$, while the scanning family $\mathcal{A} = \{A_z : z \in \mathbb{R}\}$ is just one of the examples mentioned in Section 6: $A_z = \{x \in \mathbb{R}^2 : x^1 \leq z\}$. Let $w_n$, $w_H$ be as in Section 6.

For the above regression function and the scanning family the matrix $\mathcal{C}^{-1}$ has the form

$$\mathcal{C}_z^{-1} = \frac{1}{[1-r^2][1-F(z)]} \begin{pmatrix} r^2 & -r \\ -r & 1 \end{pmatrix} + \frac{1}{1-F(z)} \begin{pmatrix} (za(z)+1)^{-1} & 0 \\ 0 & 0 \end{pmatrix}$$

and the integral in (6.3) becomes

$$\int_{-\infty}^{x_1 \wedge X_{i1}} \frac{ya(y)}{ya(y)+1} F\left(\frac{x_2 - ry}{\sqrt{1-r^2}}\right) dy\, X_{1i} - \int_{-\infty}^{x_1 \wedge X_{i1}} a(y) \frac{1}{\sqrt{1-r^2}} dy\, (X_{i2} - rX_{i1}).$$

Here, as above, $f$ and $F$ denote the standard normal density and distribution function, respectively. In our simulations the class of sets $B$ was chosen to be $(-\infty, x]$, $x \in \mathbb{R}^2$. Write $w_n(x)$, $w_H(x)$ for $w_n(B)$, $w_H(B)$ whenever $B = (-\infty, x]$, $x \in \mathbb{R}^2$, respectively. Choosing $\varphi(y) = y$ and $\tilde{\theta}$ to be the least squares estimator, the transformed process (6.3) becomes

$$w_n(x) = n^{-1/2} \sum_{i=1}^{n} \left[ \mathbb{I}(X_i \leq x) - \int_{-\infty}^{x_1 \wedge X_{i1}} \frac{ya(y)}{ya(y)+1} F\left(\frac{x_2 - ry}{\sqrt{1-r^2}}\right) dy\, X_{1i} \right.$$

$$\left. - \int_{-\infty}^{x_1 \wedge X_{i1}} a(y) \frac{1}{\sqrt{1-r^2}} dy\, (X_{i2} - rX_{i1}) \right]$$

$$\times (Y_i - \mu(X_i, \tilde{\theta})).$$

Let $\mathcal{V}_n := \sup_x |w_n(x)|$, $\mathcal{V}_H := \sup_x |w_H(x)|$.

In order to demonstrate how well the null distribution of $\mathcal{V}_n$ is approximated by the distribution of $\mathcal{V}_H$, we had to first understand the form of the latter distribution. We thus first obtained an approximation for the distribution of this r.v. as follows.

Let $H(x_1, x_2; r)$, $x = (x_1, x_2)' \in \mathbb{R}^2$, denote the d.f. of the bivariate normal distribution with standard marginals and correlation $r$. Let

(7.3)     $H_r(s,t) := H(F^{-1}(s), F^{-1}(t); r)$,     $0 \leq s, t \leq 1$,

be the corresponding copula function, and let $w(s,t) := w_H(F^{-1}(s), F^{-1}(t))$. The d.f. $P(\sup_{0 \leq s,t \leq 1} |w(s,t)| \leq v)$ is the limit as $n$ tends to infinity of, and is approximated by,

$$L_r(v) := P\left( \sup_{0 \leq s,t \leq 1} |\xi_{nH_r}(s,t) - nH_r(s,t)|/\sqrt{n} \leq v \right),$$

where $\xi_{nH_r}(s,t)$ is a Poisson process on $[0,1]^2$ with expected value $nH_r(s,t)$. Table 2 gives the simulated values of these probabilities for $r = -0.5$, $0$, $0.5$

TABLE 2
*Selected values of $(v, L_r(v))$*

(a) *for $r = -0.5$*

| $x$ | **0.71** | **0.88** | **1.00** | **1.25** | **1.50** | **1.75** | **2.00** | **2.25** | **2.50** | **2.75** | **3.00** | **3.25** |
|---|---|---|---|---|---|---|---|---|---|---|---|---|
| $L_r(x)$ | 0.00 | 0.01 | 0.05 | 0.25 | 0.50 | 0.69 | 0.82 | 0.91 | 0.95 | 0.98 | 0.99 | 0.995 |

(b) *for $r = 0$*

| $x$ | **0.66** | **0.84** | **1.00** | **1.25** | **1.50** | **1.75** | **2.00** | **2.25** | **2.50** | **2.75** | **3.00** | **3.25** |
|---|---|---|---|---|---|---|---|---|---|---|---|---|
| $L_r(x)$ | 0.00 | 0.01 | 0.07 | 0.30 | 0.53 | 0.72 | 0.84 | 0.91 | 0.95 | 0.98 | 0.99 | 0.995 |

(c) *for $r = 0.5$*

| $x$ | **0.59** | **0.79** | **1.00** | **1.25** | **1.50** | **1.75** | **2.00** | **2.25** | **2.50** | **2.75** | **3.00** | **3.25** |
|---|---|---|---|---|---|---|---|---|---|---|---|---|
| $L_r(x)$ | 0.00 | 0.01 | 0.11 | 0.35 | 0.57 | 0.74 | 0.85 | 0.92 | 0.96 | 0.98 | 0.99 | 0.996 |

and $n = 5K$, with $m = 20K$ replications. This table is based on the tables and graphs of the distribution function $L_r$ and percentile points, prepared by Dr. R. Brownrigg, available at www.mcs.vuw.ac.nz/~ray/Brownian.

Although the distribution of $\mathcal{V}_H$ depends on the copula function $H_r$, the first useful observation is that relatively sharp changes in $H_r$ do not appear to change the distribution of this r.v. by much. Table 3 summarizes a few selected percentiles to readily assess the effect of $r$ on them. It contains the values of $v_\alpha$ defined by the relation $1 - L_r(v_\alpha) = \alpha$. One readily sees that these values are very stable across the three chosen values of $r$, especially for $\alpha \leq 0.1$.

We illustrate the closeness of the distribution of $\mathcal{V}_n$ for finite $n$ to the limiting distribution with the graphs of e.d.f.s for $n = 40, 100$, with $m = 10K$ replications. Figure 2 shows the (simulated) d.f.'s of $\mathcal{V}_n$ for $n = 40, 100$; $m = 10K$, and the approximating d.f. $L_r$ (solid line) for $H_r$ as in (7.3) with $r = -0.5$, 0, and 0.5. One readily notes the remarkable closeness of these d.f.'s, especially in the right tail.

Table 4 gives the simulated values of $P(\mathcal{V}_n > v_\alpha)$ for several values of $\alpha$ and sample sizes $n = 40$ and $n = 100$, based on $m = 10K$ replications. From this table one also sees that the large sample approximation is reasonably good for even the sample size of 40 and fairly stable across the chosen values of $r$.

TABLE 3
*Selected values of $v_\alpha$ for $r = -0.5,\ 0,\ 0.5$*

| $r \setminus \alpha$ | **0.5** | **0.25** | **0.20** | **0.10** | **0.05** | **0.025** | **0.01** |
|---|---|---|---|---|---|---|---|
| $-0.5$ | 1.50 | 1.86 | 1.95 | 2.23 | 2.50 | 2.74 | 3.03 |
| $0.0$ | 1.46 | 1.81 | 1.91 | 2.21 | 2.46 | 2.70 | 3.03 |
| $0.5$ | 1.42 | 1.77 | 1.88 | 2.17 | 2.43 | 2.70 | 2.98 |

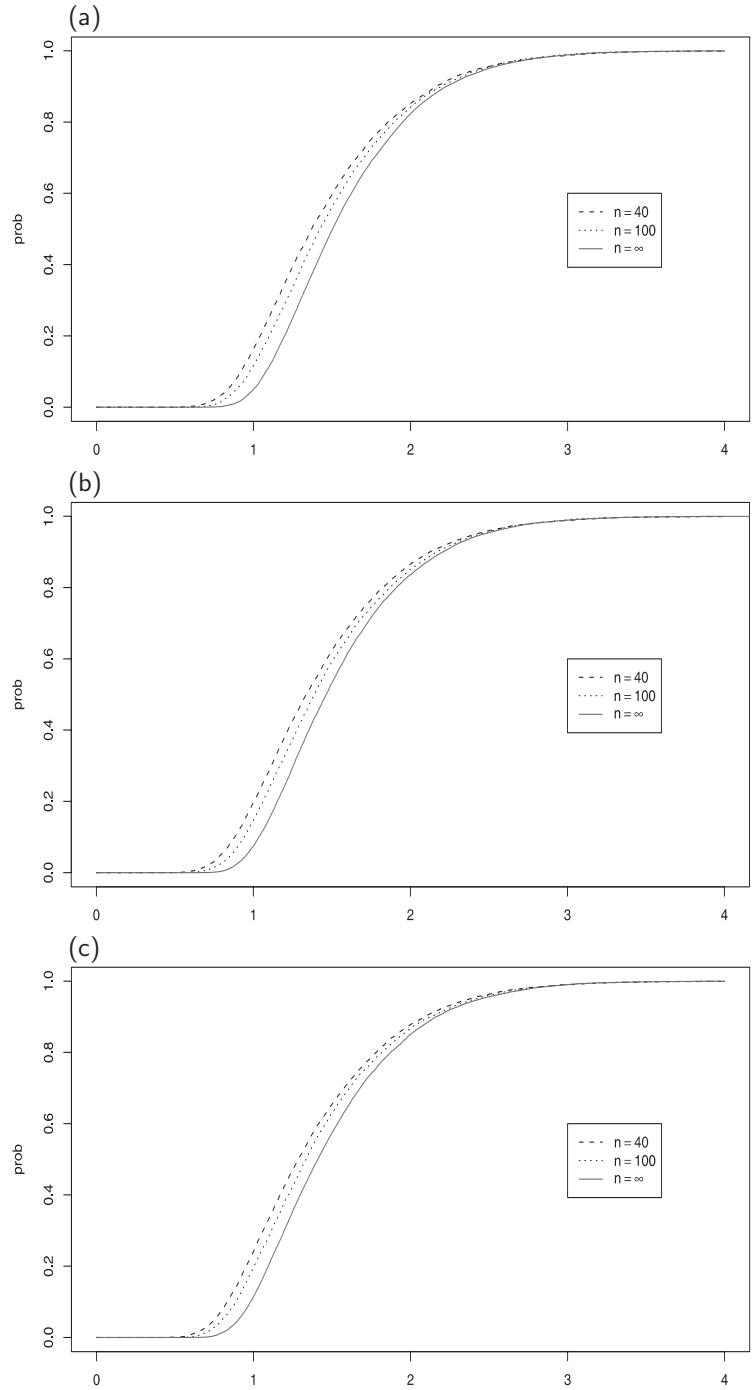

Fig. 2. (a) *E.d.f. of $\mathcal{V}_n$, $n = 40$, $100$, with $m = 10K$, and d.f. $L_r$, $r = -0.5$.* (b) *E.d.f. of $\mathcal{V}_n$, $n = 40$, $100$, with $m = 10K$, and d.f. $L_r$, $r = 0$.* (c) *E.d.f. of $\mathcal{V}_n$, $n = 40$, $100$, with $m = 10K$, and d.f. $L_r$, $r = 0.5$.*

TABLE 4
$P(\mathcal{V}_n > v_\alpha)$, $m = 20K$

| $n$ | $r \setminus \alpha$ | 0.2 | 0.1 | 0.05 | 0.01 |
|-----|------|-------|-------|-------|-------|
| 40 | $-0.5$ | 0.166 | 0.084 | 0.045 | 0.012 |
|    | 0.0 | 0.166 | 0.085 | 0.045 | 0.011 |
|    | 0.5 | 0.162 | 0.084 | 0.042 | 0.008 |
| 100 | $-0.5$ | 0.179 | 0.092 | 0.046 | 0.009 |
|     | 0.0 | 0.183 | 0.092 | 0.048 | 0.008 |
|     | 0.5 | 0.178 | 0.093 | 0.046 | 0.009 |

**Acknowledgments.** The authors thank Dr. Ray Brownrigg for his invaluable help with the numerical calculations in Section 7, and the two referees and the Associate Editor for their constructive comments.

School of Mathematical
and Computing Sciences
Victoria University of Wellington
PO Box 600
Wellington, New Zealand
e-mail: estate@mcs.vuw.ac.nz

Department of Statistics
and Probability
Michigan State University
East Lansing, Michigan 48824-1027
USA
e-mail: koul@stt.msu.edu